\documentclass[12pt,  leqno]{article}
\usepackage[english]{babel}
\usepackage[latin1]{inputenc}
\usepackage{amsmath}
\usepackage{amsfonts}
\usepackage{amssymb}
\usepackage{a4wide}

\def\Netoile{\mathbb{N} ^\ast}

\def\Qetoile{\mathbb{Q} ^\ast}

\def\trace{{\rm Trace}}

\newcommand{\Id}{{\rm Id}}

\newcommand{\N}{\mathbb{N}}

\newcommand{\Q}{\mathbb{Q}}

\newcommand{\R}{\mathbb{R}}

\newcommand{\Z}{\mathbb{Z}}

\newcommand{\al}{\alpha}

\newcommand{\capa}{\kappa}

\newcommand{\congru}{\equiv}

\newcommand{\eps}{\varepsilon}

\newcommand{\findem}{\end{pf}}

\newcommand{\moins}{\setminus}

\newcommand{\vareta}{\eta}

\newcommand{\calw}{{\cal W}}
\newcommand{\calsim}{\calr '}

\newcommand{\vect}{{\rm Span}}
\newcommand{\odu}{o(1)}
\newcommand{\oduu}{o(u_n)}

\newcommand{\zerosk}{\{0, \ldots, s_k\}}
\newcommand{\zeroskmu}{\{0, \ldots, s_k-1\}}
\newcommand{\unsk}{\{1, \ldots, s_k\}}

\newcommand{\calf}{{\cal F}}   
  
\newcommand{\calr}{{\cal R}}
\newcommand{\alphb}{{\cal A}}

\newcommand{\teta}{\vartheta}
\newcommand{\cti}{c}  

\newcommand{\Inmu}{I_{n-1}}  	
\newcommand{\uti}{\tilde u}  	
\newcommand{\deltati}{\tilde \delta}	
\newcommand{\epsti}{\tilde \eps}	

\newcommand{\beu}{\beta_1(\xi)}
\newcommand{\bez}{\beta_{0}(\xi)}
\newcommand{\beeps}{\beta_{\eps}(\xi)}

\newcommand{\bezprefpal}{\beta_{0}(\xi_{w, \varphi})}
\newcommand{\beuprefpal}{\beta_{1}(\xi_{w, \varphi})}

\newcommand{\bezpr}{\beta_0(\xi')}

\newcommand{\spcombi}{S}
\newcommand{\pal}{\pi}

\newcommand{\nbor}{\gamma}

\newcommand{\planV}{{\cal V}}
\newcommand{\hdeplanV}{H(\planV)}
\newcommand{\sevV}{{\cal V}}
\newcommand{\vk}{\sevV_k}
\newcommand{\vkmu}{\sevV_{k-1}}

\newcommand{\hdevk}{H(\vk)}
\newcommand{\hdevkmu}{H(\vkmu)}

\newcommand{\as}{\underline{a}} 
\newcommand{\ds}{\underline{d}}	
\newcommand{\es}{\underline{e}}	
\newcommand{\ws}{\underline{w}}	 
\newcommand{\vs}{\underline{v}}
\newcommand{\xs}{\underline{x}}	
\newcommand{\ys}{\underline{y}}	
\newcommand{\zs}{\underline{z}}	
\newcommand{\us}{\underline{u}} 

\newcommand{\ns}{\underline{n}^{(k)}}	 
\newcommand{\gdn}{N^{(k)}}

\newcommand{\gda}{A}

\newcommand{\gdd}{D}

\newcommand{\gdu}{U}
\newcommand{\gdx}{X}
\newcommand{\gdy}{Y}
\newcommand{\gdz}{Z}

\newcommand{\equ}{\approx}  

\newcommand{\inter}{\cap}
\newcommand{\norme}[1]{\parallel #1 \parallel }

\newcommand{\ensa}{{\cal A}}

\newcommand{\ids}{I_{0}}

\newcommand{\miroir}[1]{\widetilde{#1}}

\newcommand{\zuouvz}{(0,1]}

\newcommand{\eneq}{\end{equation}}
\newcommand{\begineq}{\begin{equation}}

\newcommand{\densi}{\delta}

\newcommand{\matri}[4]{{\tiny \left[
\begin{array}{cc} #1 & #2 \\ #3 & #4 \end{array}
\right] }}

\newtheorem{Th}{Theorem}[section]
\newtheorem{Prop}[Th]{Proposition}
\newtheorem{Cor}[Th]{Corollary}
\newtheorem{Lemme}[Th]{Lemma}

\newtheorem{Defith}[Th]{Definition} 

\newtheorem{Exemple}[Th]{\it Example\/}
\newtheorem{Remarque}[Th]{\it Remark}
\newtheorem{Step}{Step}

\def\Dem{\hspace{-\parindent}\textsc{Proof }}
\def\Demdeuxpoints{\hspace{-\parindent}\textsc{Proof: }}

\def\Question{\hspace{-\parindent}{\bf Question: }}

\newcommand{\beepsde}{\beta_{\eps_2}(\xi)}
\newcommand{\bs}{\underline{b}}
\newcommand{\espa}{\, \, \, \, ;  \hspace{0.5cm}}
\newcommand{\espabis}{\hspace{0.3cm}}

\newcommand{\Span}{{\rm Span}}

\title{Palindromic Prefixes and Diophantine Approximation}
\author{St\'ephane Fischler}
\date{\today}

\begin{document}

\maketitle

{\bf Abstract:}
This text is devoted to simultaneous approximation to $\xi$ and $\xi^2$ by rational numbers with the same denominator, where $\xi$ is an irrational  non-quadratic real number. We focus on an exponent $\bez$ that measures the regularity of the sequence of all exceptionally precise  such approximants. We prove that $\bez$ takes the same set of values as a combinatorial quantity that measures the abundance of palindromic prefixes in an infinite word $w$. This allows us to give a precise exposition of Roy's palindromic prefix method. The main tools we use are Davenport-Schmidt's sequence of minimal points and Roy's bracket operation.

\bigskip

\bigskip


\bigskip

\bigskip


\section{Introduction} \label{sec1}

Throughout this text, we denote by $\xi$ a real number, assumed to be irrational  and non-quadratic (that is, $[\Q(\xi) : \Q] \geq 3$). We study the quality of simultaneous rational approximants to $\xi$ and $\xi^2$, that is the possibility to find triples $\xs = (x_0, x_1, x_2) \in \Z^3$ with
$$L(\xs) = \max ( | x_0 \xi - x_1|, |x_0 \xi^2 - x_2|)$$
very small in comparison with $|x_0|$. In more precise terms, for any $0 < \eps \leq 1$ we consider the exponent $\beeps$ defined as follows: $\beeps$ is the infimum of the set of all $\beta$ such that for any sufficiently large $B > 0$ there exists $\xs \in \Z^3$ such that 
\begin{equation} \label{eqp1}
1 \leq | x_0 | \leq B \mbox{ and } L(\xs) \leq \min (B^{-1/\beta} , |x_0|^{\eps-1}).
\eneq
If this set of $\beta$ is empty, we let $\beeps = + \infty$. For $\eps = 0$, we let 
$$\bez = \sup_{0 < \eps < 1} \beeps = \lim_{\eps \to 0} \beeps$$
since for any fixed $\xi$, the map $\eps \mapsto \beeps$ is non-increasing. We do not consider $\beeps$ for $\eps < 0$, since $L(\xs) \leq | x_0 |^{\eps - 1} $ with $\eps < 0$ implies $\frac{x_2}{x_0} = (\frac{x_1}{x_0})^2$ if $|x_0|$ is sufficiently large (see \cite{DS}); therefore these exponents merely concern rational approximants to $\xi$.

The exponents $\beeps$, introduced in \cite{CRASasim}, generalize the classical exponent $\beu$ (denoted by $\hat w'_2(\xi)$ or by $1/ \hat \lambda_2(\xi)$ in \cite{BL}, \cite{Bugeaudlivre} 
and \cite{Roytwoexp}) for which the following results are known:
\begin{itemize}
\item $\beu \leq 2$ for any $\xi$, by applying Dirichlet's pigeon-hole principle;
\item $\beu \geq \nbor = \frac{1 + \sqrt5}{2} = 1.618 \ldots$ for any $\xi$, proved by Davenport and Schmidt \cite{DS};
\item There exists $\xi$ such that $\beu = \nbor$, proved by Roy \cite{RoyCRAS};
\item The set of values taken by $\beu$ is dense in $[\nbor, 2]$, proved by Roy \cite{Roytwoexp}. 
\end{itemize}

\bigskip

The aim of this paper is to prove analogous results on $\bez$. One can prove easily (see \cite{BL}) that $\bez = + \infty$ for almost all $\xi$ with respect to Lebesgue measure, and Davenport-Schmidt's lower bound on $\beu$ yields $\bez \geq \nbor$ for any $\xi$. Moreover Roy's {\em palindromic prefix method} (\cite{RoyCRAS}; see also \cite{CRASasim}) allows one to obtain a $\xi$ such that $\bez < 2$, starting from any word $w$ with ``sufficiently many'' palindromic prefixes. In more precise terms, let $w = w_1w_2 \ldots$ be an infinite word on a finite alphabet. We denote by $(n_i)_{i \geq 1}$ the increasing sequence of all lengths of palindromic prefixes of $w$, in such a way that $n$ is among the $n_i$ if, and only if, the prefix $w_1 w_2 \ldots w_n$ of $w$ is a palindrome (i.e., $w_1 = w_n$, $w_2 = w_{n-1}$, \ldots). We let 
$$\densi(w) = \limsup_{i \to + \infty} \frac{n_{i+1}}{n_i}$$
if the sequence $(n_i)$ is infinite, and $\densi(w) = + \infty$ otherwise. The words $w$ such that $\densi(w) < 2$ are studied in \cite{SFcombi}. Roy's palindromic prefix method enables one to construct, from any non ultimately periodic word $w$ such that $\densi(w) <2$, a real number $\xi$ such that $\bez < 2$. Moreover, the equality $\bez = \densi(w)$ holds if $w$ is the Fibonacci word (see \cite{RoyCRAS}), and more generally for any characteristic Sturmian word $w$ (see \cite{BL}). In this paper, we prove this equality for any non ultimately periodic word $w$ such that $\densi(w) <2$. We also show a reciprocal statement:

\begin{Th} \label{th1}
We have 
$$
\{ \bez \, \,  ; \espabis \xi \in \R \mbox{  irrational,   non-quadratic }\} \cap (1,2) 
=
\{ \densi(w)  \, \, ; \espabis w \mbox{ infinite word}\} \cap (1,2) .$$
\end{Th}

The set $\spcombi = \{ \densi(w) \} \cap (1,2) $ that appears in Theorem \ref{th1} is studied in \cite{SFcombi}. The least element of $\spcombi$ is the golden ratio $\nbor$ (as implied by the previously mentioned results). Apart from $\nbor$, the least element is $\sigma_2 =1 + \frac{\sqrt2}{2} = 1.707\ldots$, then (apart from $\sigma_2$) it is $\sigma_3 = \frac{2 + \sqrt{10}}{3}$, and so on: there is an increasing sequence $(\sigma_n)$ of isolated points in $\spcombi$, that converges to the least accumulation point $\sigma_{\infty} = 1.721\ldots$ of $\spcombi$. This follows from a result of Cassaigne \cite{Cassaigne} and the property (\cite{SFcombi}, Theorem 1.3) that $\spcombi \inter (1, \sqrt3]$ coincides with the set of values of $\densi(w)$ coming from characteristic Sturmian words $w$. Combining this property with Theorem \ref{th1} proves Th\'eor\`eme 2.1 announced in \cite{CRASasim}, namely: all values of $\bez$ less than $\sqrt3$ can be obtained from characteristic Sturmian words as in \cite{BL}. In particular,   the following corollary shows that the situation is completely different from the case of the exponent $\beu$, which assumes \cite{Roytwoexp} a set of values dense in $[\nbor , 2]$. 

\begin{Cor} There is no real number $\xi$ such that $\nbor < \bez < \sigma_2 = 1.707\ldots$.
\end{Cor} 

\bigskip

We shall deduce Theorem \ref{th1} from a ``structure theorem'' on numbers $\xi$ such that $\bez < 2$, namely that the sequence of all ``exceptionally precise'' approximants to $\xi$ and $\xi^2$ satisfies the ``same'' recurrence relation as the sequence of all palindromic prefixes of some word $w$ with $\densi(w) < 2$. This enables us to associate with each $\xi$ such that $\bez < 2$ a word $w$ such that $\densi(w) < 2$, in such a way that $\bez = \densi(w)$. This gives a kind of converse to Roy's palindromic prefix method. As byproduct of our approach, we also obtain the following result:

\begin{Th} \label{thcor516}
Let $\xi$ be such that $\bez < 2$. Then there exists $\eps_1 > 0$ (depending only on $\bez$) such that $\beeps = \bez$ for any $\eps < \eps_1$.
\end{Th}

\bigskip

Throughout the text, we shall use the following notation. We denote by $\N = \{0, 1, 2, \ldots \}$ the set of non-negative integers, and let $\Netoile  = \N \moins \{0\}$.  We write  $u_t \ll v_t$ (and $v_t \gg u_t$) if $v_t$ is positive for $t$ sufficiently large, and $u_t/v_t$ is bounded from above as $t$ tends to infinity. When $u_t \ll v_t$ and $v_t \ll u_t$, we write $u_t \equ v_t$. At last, when a point of $\Z^3$ is referred to as an underlined letter (e.g. $\xs$), we denote by  the corresponding uppercase letter its norm, that is the largest absolute value of its  coordinates (e.g. $X = \max(|x_{0}|, |x_{1}|, |x_{2}|) $ if $\xs = (x_0, x_1, x_2)$). 

\bigskip

The structure of this text is as follows.
In Section \ref{sec2bis} we recall the notation and results of \cite{SFcombi} about words with many palindromic prefixes, and also the definition and properties of Roy's bracket operation. Our main tool is the sequence of minimal points introduced by Davenport and Schmidt. Section  \ref{sec2} is devoted to this sequence: we recall classical results and prove new ones (which may be of independent interest). We state and prove in Section \ref{sec5} our main diophantine result, which is the crucial step in the proof (\S \ref{subsec53}) of the results stated in this introduction. This enables us in Section \ref{sec6} to give a precise account on Roy's  palindromic prefix method -- which was the original motivation of this paper. At last, Section \ref{sec7} is devoted to some open questions.

\bigskip

\noindent{\bf Acknowledgements: } This work has begun during a 6-month stay at Ottawa University. I am very grateful to Damien Roy for this invitation, to Michel Waldschmidt for useful discussions, and to the referee for important improvements in the redaction.

\section{Palindromes, Roy's Bracket and the Main Results} \label{sec2bis}

In this Section, we recall the properties  \cite{SFcombi} of words with many palindromic prefixes, and also the bracket operation introduced by Roy \cite{RoyPLMS}.  

\subsection{Words with Many Palindromic Prefixes} \label{subsec26}

Let us recall the definitions and results of \cite{SFcombi} that will be useful here (with minor changes intended to fit into the diophantine setting). 

\bigskip

Let $\psi : \Netoile \to \N $ be a function such that 
\begineq \label{eq261}
\psi(n) \leq n-1 \hspace{1cm} \mbox{ for any } n \geq 1. 
\eneq
Denote by $(\teta_k)_{k \geq 0}$ the sequence  of all indices $n \geq 1$ (in increasing order) such that $\psi(n) \leq n-2$. This sequence  may be either finite or infinite; it  is infinite for the functions $\psi$ of interest here (see \eqref{eq3} below). We shall use the following variant of Definition 4.9 of \cite{SFcombi}:

\begin{Defith} \label{def262}
A function $\psi$ is said to be {\em asymptotically reduced} if \eqref{eq261} holds and the associated sequence $(\teta_k)$ is such that 
$$\psi(\teta_k)< \teta_{k-1} \mbox{ and } \psi(\teta_k) \neq \psi (\teta_{k-1})$$
for any sufficiently large $k $. When the sequence  $(\teta_k)$ is finite, we agree  that $\psi$ is asymptotically reduced.
\end{Defith}

Let us denote by $\calf$ the set of all functions $\psi : \Netoile \to \N $ such that 
\begin{equation} \label{eq3}
\left\{
\begin{array}{l}
\psi \mbox{ is asymptotically reduced} \\
\mbox{There exists $\cti$ such that } i - \cti \leq \psi(i) \leq i-1 \mbox{ for any } i \geq 1\\
\mbox{There are infinitely many $i$ such that } \psi(i) \leq i-2.
\end{array}
\right.
\end{equation}
For $\psi \in \calf$, the sequence $(\teta_k)$ is infinite, and we have $\teta_{k+1} \leq \teta_k + c-1$ for any $k$ sufficiently large (since $\teta_{k+1} - c \leq \psi(\teta_{k+1}) < \teta_k$). 

\bigskip

Given an infinite word $w = w_1w_2\ldots$ on an arbitrary (finite or infinite) alphabet $\alphb$, we denote by $(\pal_i)_{i \geq 1}$ the (finite or infinite) sequence of all palindromic prefixes of $w$, and by $n_i$ the length of $\pal_i$ (in such a way that $\pal_i = w_1 w_2 \ldots w_{n_i}$ is a palindrome for any $i \geq 1$).
For $i' \leq i$, $\pal_{i'}$ is a prefix of $\pal_i$, so there is a word $b = w_{n_{i'}+1} \ldots w_{n_i}$ such that $\pal_i = \pal_{i'} b$. We denote by $\pal_{i'}^{-1} \pal_i$ this word $b$. 

\smallskip

 We let $\calw$ be the set of all non ultimately periodic words $w$ such that the increasing sequence $(n_i)$ is infinite and satisfies $\limsup n_{i+1}/n_i < 2$. Then the following result is proved in \cite{SFcombi}:

\begin{Th} \label{thcombi}

$(i)$ For any $w \in \calw$ there exists $\psi \in \calf$ such that
\begineq \label{eq2}
\pal_{i+1} = \pal_i \pal_{\psi(i)}^{-1} \pal_i \mbox{ for any $i$ sufficiently large.}
\eneq

$(ii)$ For any $\psi \in \calf$ there exists $w \in \calw$ such that \eqref{eq2} holds.
\end{Th}

Theorem \ref{thcombi} implies that any $w \in \calw$ can be written on a finite alphabet. So we may assume, throughout this text, that the alphabet $\alphb$ is finite.

The following relation follows from \eqref{eq2} and will be used repeatedly:
$$n_{i+1}  = 2 n_i - n_{\psi(i)}  \mbox{ for any $i$ sufficiently large.}$$
 
\medskip

The three conditions that appear in \eqref{eq3} are of different nature. The last one corresponds, in the definition of $\calw$, to the assumption that $w$ is not ultimately periodic (which is equivalent to $w$ not being periodic: see Lemma 5.6 of \cite{SFcombi}). The second one enables us to get rid of the case $\densi(\psi) = 2$; this is very useful since  Theorem \ref{thcombi} generalizes only to a specific class of words $w$ such that $\densi(w) = 2$. At last, the assumption that $\psi$ is asymptotically reduced ensures that distinct functions $\psi$ (modulo the equivalence relation $\calr$ defined below) always correspond to distinct words $w$.

\begin{Exemple} \label{explefibo}
The Fibonacci word $abaaba \ldots$ on the two-letter alphabet $\{a,b\}$ corresponds (in Theorem \ref{thcombi}) to all functions $\psi$ such that $\psi(n) = n-2$ for any sufficiently large $n$ (see \cite{SFcombi}, Example 4.7). More generally, if $w$ is the characteristic Sturmian word with slope $[0, s_1, s_2, \ldots]$ then $\densi(w) = 2$ if, and only if, the sequence $(s_n)$ is unbounded. In the opposite case, $w$ belongs to $\calw$ and corresponds to any function $\psi$ such that $\psi(n) = n-s_k - 1$ if $n$ can be written $s_1 + \ldots + s_k$ for some $k$, and $\psi(n) = n-1$ otherwise, for $n$ sufficiently large (see \cite{SFcombi}, Example 4.6). We recover the Fibonacci word by considering the special case $s_1 = s_2 = \ldots = 1$.
\end{Exemple}

\begin{Exemple} \label{expleto}
Let $w \in \calw$ be written on a finite alphabet, say $\alphb = \{a_1, \ldots, a_r\}$. Let $p_1, \ldots, p_r$ be arbitrary palindromes, written on another alphabet $\alphb'$. Then replacing each $a_i$ with $p_i$ in the word $w$ yields a word $w'$ written on $\alphb'$. Denote by $(\pal_i)_{i \geq 1}$ the sequence of all palindromic prefixes of $w$, and by $\pal'_i$ the finite word obtained from $\pal_i$ by replacing each $a_i$ with $p_i$. Since $p_1, \ldots, p_r$ are palindromes, all $\pal'_i$ are palindromic prefixes of $w'$. However, in general, the sequence $(\pal'_i)_{i \geq 1}$ does not contain {\em all} palindromic prefixes of $w'$. For instance, if all $p_i$ have length greater than 2 then the first letter of $w'$ is a palindromic prefix of $w'$ but is not among the $\pal'_i$. 
\end{Exemple}

Throughout this text, we are interested only in asymptotic properties, so the palindromic prefixes of $w'$ that are missing in Example \ref{expleto} are not a problem (as long as their number is finite). To take this observation into account, we define an equivalence relation $\calr$ on $\calf$ as follows: 

\begin{Defith} \label{deficalr}
For $\psi, \psi' \in \calf$, we set $\psi \calr \psi'$ if there exist $\delta$ and $i_1$ such that $\psi(i) - i = \psi'(i-\delta) - (i - \delta)$ for any $i \geq i_1$.
\end{Defith}

We now define an analogous equivalence relation on $\calw$:

\begin{Defith} \label{deficalsim}
Let $w, w' \in \calw$, and $(\pal_i)_{i \geq 1}$ (resp. $(\pal'_i)_{i \geq 1}$) be the sequence of all palindromic prefixes of $w$ (resp. $w'$). We set $w \calsim w'$ if there exist $\psi \in \calf$ and $\delta$ such that 
$$\pal_{i+1} = \pal_i \pal_{\psi(i)}^{-1} \pal_i 
\mbox{ and }
\pal'_{i+1} = \pal'_i {\pal'}_{\psi(i-\delta)+\delta}^{-1} \pal'_i 
$$
for any $i$ sufficiently large.
\end{Defith}

This definition means that $w \calsim w'$ if, and only if, there exist $\psi, \psi' \in \calf$ with $\psi \calr \psi'$ such that 
$\pal_{i+1} = \pal_i \pal_{\psi(i)}^{-1} \pal_i $
and
$\pal'_{i+1} = \pal'_i {\pal'}_{\psi'(i)}^{-1} \pal'_i $
for any $i$ sufficiently large. Another way of stating this is the following: $w \calsim w'$ if, and only if, after omitting a finite number of initial terms in the sequence $(\pal_i)_{i \geq 1}$ or in  $(\pal'_i)_{i \geq 1}$, we have $\pal_{i+1} = \pal_i \pal_{\psi(i)}^{-1} \pal_i $
and
$\pal'_{i+1} = \pal'_i {\pal'}_{\psi(i)}^{-1} \pal'_i $ for some $\psi \in \calf $ and any 
$i$ sufficiently large.

\bigskip

Thanks to Definitions \ref{deficalr} and \ref{deficalsim}, Theorem \ref{thcombi} yields a bijective map
\begin{equation} \label{eqbij}
\calw / \calsim \stackrel{\sim}{\longrightarrow} \calf / \calr. 
\end{equation}

\bigskip

We now define a quantity $\densi(\psi)$ as follows:
\begin{Defith} 
For $\psi \in \calf$, we let 
$$\densi(\psi) = \limsup \frac{m_{i+1}}{m_i}$$
where $(m_i)_{i \geq 1}$ is any increasing sequence of non-negative integers such that $m_{i+1} = 2m_i - m_{\psi(i)}$ for any $i$ sufficiently large.
\end{Defith}

The value of $\densi(\psi)$ does not depend on the choice of a peculiar sequence $(m_i)$, thanks to Proposition 6.2 of \cite{SFcombi}. In particular, if $\psi$ corresponds to a word $w$ as in Theorem \ref{thcombi}  then one can choose $m_i = n_i$, hence $\densi(\psi) = \densi(w)$. 
It follows from this remark and from Proposition 6.2 of \cite{SFcombi} that :
\begin{itemize}
\item  The map $\calw \to (1,2), \, w \mapsto \densi(w)$ factors into a map $\calw / \calsim \to (1,2)$.
\item  The map $\calf \to (1,2), \, \psi \mapsto \densi(\psi)$ factors into a map $\calf / \calr \to (1,2)$.
\item These two maps, together with the bijection \eqref{eqbij}, make up a commutative diagram:
\begin{equation} \label{eqdiag}
\begin{array}{rcccl}
\calw / \calsim &  & \stackrel{\sim}{\longrightarrow} &   & \calf / \calr \\
			& \searrow& 		& \swarrow&  \\
			& & (1,2) & 		&  
\end{array}
\end{equation}
\end{itemize}

We have constructed in \cite{SFcombi}  (Remark 7.6) two words $w$ and $w'$ such that $\densi(w) = \densi(w') = \sqrt3$, but with $w \not\congru w' \mod \calsim$ (actually $w$ is not episturmian whereas $w'$ is characteristic Sturmian). This proves that in the diagram \eqref{eqdiag}, 
the maps  $\calw / \calsim \to (1,2)$ and $\calf  / \calr \to (1,2)$  induced by $\densi$  are not injective.

\subsection{Classical Estimates about Approximants} \label{subsec22}

Throughout this text, we let for any $\eps$ such that  $0 < \eps < 1$:
$$\ensa_\eps = \{\xs \in \Z^3  \moins \{0\}  \espa  L(\xs) \leq \gdx^{-(1-\eps)}\}, $$
with $X = \max(|x_{0}|, |x_{1}|, |x_{2}|) $ (as explained in the Introduction). 
By convention, for $\eps \geq 1$,  we let $\ensa_{\eps} = \Z^3 \moins \{0\}$. Moreover, 
we identify a point $\xs = (x_0,x_1,x_2) \in \Z^3$ with the symmetric matrix
${\tiny \left[ \begin{array}{cc} x_0 & x_1 \\ 
x_1 & x_2 \end{array} \right]}$. 

\bigskip

In this paragraph, we state classical estimates which  are  due (mainly) to Davenport and Schmidt.
For instance, the first inequality of  the following lemma is proved in Lemmas 2 and 3 of \cite{DS}, the second one in Lemma 4.

\begin{Lemme} \label{lem23}
We have $$\vert \det(\xs) \vert = | x_0 x_2 - x_1^2| \ll \gdx L(\xs).$$
Moreover, if $\xs, \ys, \zs$ are such that $\gdx \geq \max(\gdy,\gdz)$ and $L(\xs) \leq \min(L(\ys),L(\zs))$, then we have
$$\vert \det(\xs, \ys, \zs) \vert \ll \gdx L(\ys) L(\zs). $$
\end{Lemme}

\bigskip

Let $\norme{\xs}$ denote the norm of a vector $\xs \in \Z^3$ (that is, its greatest coordinate in absolute value, usually denoted by $X$ in this text). 
Let $\planV$ be a sub-$\Z$-module of $\Z^3$, of rank 2. Then
$\norme{\xs \wedge \ys}$, computed for a $\Z$-basis $(\xs,\ys)$ of $\planV$, 
does not depend on the chosen basis but only on $\planV$. This is the {\em height} of $\planV$, denoted by $\hdeplanV$. 

Let $\xs$ and $\ys$ be two linearly independent vectors in 
$\planV$. Then
$\norme{\xs \wedge \ys} = N \, \hdeplanV$ where $N$ is the index, in $\planV$, of the subgroup generated by   $\xs$ and $\ys$. Moreover, Lemma 3 of  \cite{DS} shows that 
\begineq \label{eq24} 
 \norme{\xs \wedge \ys} \ll X L(\ys) + Y L(\xs) .
\eneq

\subsection{Roy's Bracket Operation} \label{subsec21}

The bracket operation $[.,.,.]$ introduced by Roy \cite{RoyPLMS} will be used in a crucial manner in the statement, and proof, of our results. If 
$\xs,\ys,\zs$ are three linearly dependent vectors in $\Z^3$, understood as symmetric matrices, then the matrix $-\xs J \zs J \ys$ (where $J={\tiny \left[ \begin{array}{cc} 0 & 1 \\ -1 & 
0 \end{array} \right]}$) is also symmetric, and it is denoted by 
$[\xs, \ys, \zs]$. The importance of this operation in our context is easily seen, for instance, by considering Lemma \ref{lem64} below. The following relation holds:
$$\det([\xs, \ys, \zs]) = \det(\xs) \, \det(\ys) \, \det(\zs). $$
Moreover $\zs$ can be obtained back from  $\xs$, $\ys$ and $[\xs, \ys, \zs]$ using the following formula:
\begineq \label{eq211}
\det(\xs) \det(\ys) \zs = [\xs,\ys,[\xs,\ys,\zs]]
\eneq
which makes sense since $\xs$, $\ys$ and$[\xs,\ys,\zs]$ are linearly dependent (see \cite{RoyPLMS}, Lemma 2.1). At last, since $(J \ys)^2 = - \det (\ys) \Id$ we have:
\begineq \label{eq212}
 [\xs,\ys, \ys ] = \det(\ys) \xs.
\eneq

\medskip

The main property of this bracket is the following 
 slight generalization of \cite{RoyPLMS}, Lemma 3.1 $(iii)$ (with $X = \max(|x_{0}|, |x_{1}|, |x_{2}|)$ and so on).

\begin{Prop} \label{prop21}
Let 
$\xs,\ys,\zs$ be linearly dependent vectors in $\Z^3$.
Let 
$$\lambda = \gdz L(\xs) L(\ys) + L(\zs) \min(\gdy L(\xs), \gdx L(\ys)). $$
Then the bracket $\us = [\xs, \ys , \zs]$ is such that 
$$\gdu \ll \gdx \gdy L(\zs) + \lambda
\mbox{ and }
L(\us) \ll \lambda. $$
\end{Prop}

Proposition \ref{prop21} will also be used through the following Corollary:

\begin{Cor} \label{cor22}
Let $\capa, \eps, \eps' \in (0,1)$ be such that $\eps \frac{1 + \capa}{1 - \capa} < \eps' < 1$.
Let $\xs,\ys,\zs$ be linearly dependent vectors in 
$\ensa_\eps$, with $X$, $Y$, $Z$ large enough in terms of $\capa$, $ \eps$, $\eps'$ and 
$\xi$. We assume 
$$\gdx, \gdy \leq \gdz \mbox{ , } \hspace{1cm}
\gdz \leq (\gdx \gdy)^\capa \hspace{0.6cm} \mbox{ and } \hspace{0.6cm}
\det(\xs), \det(\ys), \det(\zs) \neq 0. $$
Then the bracket  $\us = [\xs, \ys , \zs]$ belongs to $\ensa_{\eps'}$.
\end{Cor}

\Demdeuxpoints  
Let the notation be as in the Proposition \ref{prop21}. Then we have $\lambda \ll Z (XY)^{\eps-1} \ll (XY)^{\capa + \eps-1} <  1$ since $\eps+ \capa < 1$, and $X Y L(\zs) \gg XY Z^{-1} \geq (XY)^{1-\capa} >  1$ since $\det(\zs) \neq 0$ (using Lemma \ref{lem23}). Therefore $U \ll X Y L(\zs) $ and $L(\us) \ll \lambda$, hence
$$U^{1 - \eps'} L(\us) \ll (XY)^{\eps - \eps'} Z^{\eps + \eps' - \eps \eps'} \leq (XY)^{(1+\capa)\eps - (1-\capa)\eps'}.$$
Since $\xs$, $\ys$ and $\zs$ have non-zero determinants, $\us$ is not zero; this concludes the proof of Corollary \ref{cor22}.

\section{Davenport-Schmidt's Minimal Points}  \label{sec2}

In this Section, we recall (\S \ref{subsec23})  the definition and classical properties of the sequence of minimal points, introduced by Davenport and Schmidt \cite{DS}. Then we move to new results on this sequence (\S\S \ref{subsec32} to \ref{subsec41}),  that are of independent interest and may be used to study numbers $\xi$ such that $\beu  <2$  but not necessarily $\bez < 2$. The most important results are  Proposition \ref{prop36} (which states that any sufficiently precise approximation is collinear to a minimal point) and Proposition \ref{prop41nv} (which deals with linear independence of approximants).

\bigskip

We shall make a frequent use of the notation and results of \S \ref{subsec22}. 
Recall that for any $\eps$ such that  $0 < \eps < 1$, we let 
$$\ensa_\eps = \{\xs \in \Z^3  \moins \{0\}  \espa  L(\xs) \leq \gdx^{-(1-\eps)}\}, $$
with $X = \max(|x_{0}|, |x_{1}|, |x_{2}|) $, and
 $\ensa_{\eps} = \Z^3 \moins \{0\}$   for $\eps \geq 1$.  Moreover, 
we identify a point $\xs = (x_0,x_1,x_2) \in \Z^3$ with the symmetric matrix
${\tiny \left[ \begin{array}{cc} x_0 & x_1 \\ 
x_1 & x_2 \end{array} \right]}$. 

\bigskip

Throughout this Section, we let $\xi$ denote an irrational  non-quadratic real number such that $\beu < 2$. We do not assume anything about $\bez$.

\subsection{Definition and Notation} \label{subsec23}

Let $\xi$ be an irrational  non-quadratic real number. For any real $X \geq 1$, the set of all
$\as = (a_0,a_1,a_2) \in \Z^3$ such that 
 $1 \leq a_0 \leq X$ and $L(\as) \leq 1$ is finite, and contains exactly one element for which
 $L$ is minimal. Following \cite{DS67} and \cite{DS}, we call this element 
 {\em minimal point  corresponding to $X$}. We denote by 
 $(\as_i)_{i \geq 1}$ the sequence of all minimal points, in such a way that 
 $\as_i$ be a minimal point  corresponding to all $X$ such that 
$a_{i,0} \leq X < a_{i+1,0}$. By definition, any $\as_i$ is a primitive point of 
$\Z^3$ (that is,  $a_{i,0}$, $a_{i,1}$ and $a_{i,2}$ are globally coprime), and any two distinct minimal points are always linearly independent over $\Z$. As usual, we let $A_i =  \max(|a_{i,0}|, |a_{i,1}|, |a_{i,2}|) \geq 1$. For any $i$ sufficiently large, we have $L(\as_i) < 1/2$ so that $a_{i,1}$ (resp. $a_{i,2}$) is the closest integer to $a_{i,0} \xi$ (resp. $a_{i,0} \xi^2$), and $A_i < A_{i+1}$. 

\medskip

From now on, we assume $\beu < 2$. Then Lemma 2 of \cite{DS} implies 
$\det(\as_i) \neq 0$ for $i$ sufficiently large, hence
\begineq \label{eq231} 
1 \ll \gda_i L(\as_i)
\eneq
using Lemma  \ref{lem23}.

Let $\ids$ denote the set of all indices $i$ such that  $\as_{i-1}$, $\as_i$ and
$\as_{i+1}$ are linearly independent. This set is infinite (see \cite{DS}); we denote by $(i_k)$ the sequence of all elements of $\ids$, in increasing order. We let 
$$\ds_k = \as_{i_k}. $$

\smallskip

We denote by $\vk$ the intersection with $\Z^3$ of the sub-$\Q$-vector space of $\Q^3$ spanned by  $\ds_k$ et $\ds_{k+1}$, and by $\hdevk$ its height (defined in \S \ref{subsec22}). It follows from \cite{DS} that $\hdevk$ tends to infinity as $k$ tends to infinity. Moreover, Lemma 4.1 of \cite{RoyPLMS} (see also Lemma 2 of \cite{DS67}) is the following result:

\begin{Lemme} \label{lem3131}
For any $k$ and $i$ sufficiently large such that
$i_{k} \leq i < i_{k+1}$, the points $\as_i$ and $\as_{i+1}$ make up a basis of the 
$\Z$-module $\vk$ and we have $\gda_{i+1} L(\as_i) \equ \hdevk$.
\end{Lemme}

Let $\eps \in \zuouvz$, and $(\us_i)_{i \geq 1} $ be the 
 sequence of all minimal points in $\ensa_{\eps}$, ordered according to their norms $U_i$. Then, we have 
\begin{equation} \label{eqscan}
\beeps = \limsup_{i \to \infty} \frac{\log U_{i+1}}{- \log L(\us_i)}.
\end{equation}
This fact immediately follows from the definition of $\beeps$ (see \cite{DS}); the point is that in Equation \eqref{eqp1} we may assume that $\xs$ is a minimal point.

\subsection{Estimates for the Height of $\vk$} \label{subsec32}

In the following estimates, we assume $\beu < 2$ and use repeatedly the following consequence of Equation \eqref{eqscan}: for any $\al > \beu $ we have 
\begineq \label{eq321} 
\gda_{i+1} \ll L(\as_i) ^{-\al} \mbox{ and }
L(\as_i) \ll \gda_{i+1} ^{-1/\al}. 
\eneq


\begin{Lemme} \label{lem32}
Let $\eps \in \zuouvz$ and 
$i$ be such that $a_i \in \ensa_\eps$. Then we have 
$$\gda_{i+1} \gg \hdevk \gda_i ^{1-\eps}$$
where $k$ is the only integer such that $i_k \leq i < i_{k+1}$.
\end{Lemme}

\Demdeuxpoints We have $\hdevk \equ \gda_{i+1} L(\as_i) \leq \gda_{i+1} \gda_i ^{-1+\eps}$, thanks to Lemma \ref{lem3131} and by definition of $ \ensa_\eps$. 

\bigskip

\begin{Lemme} \label{lem33}
For any $\al > \beu $ we have 
$$\hdevk \gg \gdd_{k+1} ^{2 - \al}. $$
\end{Lemme}

\Demdeuxpoints 
Let $i = i_{k+1}$; then the following inequalities hold thanks to Lemmas \ref{lem23} and \ref{lem3131}:
\begin{eqnarray*}
1 &\ll & \vert \det (\as_{i-1}, \as_i, \as_{i+1}) \vert \\
&\ll & \gda_{i+1} L(\as_i) L(\as_{i-1}) \\
&\ll & L(\as_i)^{1-\al} \gda_i ^{-1} \hdevk \\
&\ll & \gda_i ^{\al-2} \hdevk .
\end{eqnarray*}
This proves the Lemma.

\bigskip

\begin{Lemme} \label{lem34}
For any $\al > \beu $ we have
$$\hdevk \gg \gdd_k ^{1/\al}. $$
\end{Lemme}

\Demdeuxpoints 
Let $i = i_k$. The following inequalities hold:
\begin{eqnarray*}
1 &\ll & \vert \det (\as_{i-1}, \as_i, \as_{i+1}) \vert \\
&\ll & \gda_{i+1} L(\as_i) L(\as_{i-1}) \\
&\ll & \hdevk L(\as_{i-1}) \\
&\ll & \hdevk \gda_i ^{-1/\al} 
\end{eqnarray*}
and prove the Lemma. 

\bigskip

\begin{Remarque} \label{rem34bis}
Lemmas \ref{lem33} and \ref{lem34} are optimal when $\xi$ is the number constructed by Roy \cite{RoyCRAS} from the Fibonacci word. Actually, for this number we have 
$\hdevk \equ \gdd_{k+1} ^{2 - \nbor} \equ \gdd_k ^{1/\nbor}$. Now, let us assume that $\xi$ is, more generally, constructed (as in \cite{Alloucheetal}) and  \cite{BL}) from a characteristic Sturmian word $w$. Then Lemma \ref{lem33} is still optimal for some values of $k$; actually, if $k$ is such that $A_{i_{k+1}+1} = L(\ds_{k+1})^{-\al}$ then the estimates proved by Bugeaud and Laurent imply $\hdevk \equ D_{k+1}^{2-\al}$. However, Lemma \ref{lem34} is not optimal, as can already be seen from the proof: the upper bound $L(\as_{i_k-1}) \ll D_k ^{-1/\al}$ can be sharp only if $\as_{i_k-1} = \ds_{k-1}$, and there are words $w$ for which this never happens. 
\end{Remarque}

\subsection{Importance of Minimal Points} \label{subsec33}

By definition of $\beu $, for any $\eps > 1 - 1/\beu $ all minimal points belong to $\ensa_\eps$ (except for a finite number of exceptions, depending on $\eps$). This means that any 
minimal point is a rather good simultaneous approximant to  $ \xi$ and $\xi^2$. The following Proposition is a converse statement, and proves that the approximations provided by the palindromic prefix method (see \S \ref{subsec62}) are minimal points (up to proportionality).

\begin{Prop} \label{prop36}
Let $\eps > 0$ be such that $\eps < 2 - \beu $.
Let $\xs \in \ensa_\eps$ with $\gdx$ sufficiently large (in terms of $\eps$ and $\xi$). Then $\xs$ is collinear to some minimal point $\as_i$. 
\end{Prop}

\Question Is it possible, in Proposition \ref{prop36}, to weaken the assumption on $\eps$ to $\eps < 1 - 1/\beu$ ? 

\bigskip

A positive answer to this question would be very satisfactory, since the assumption on $\eps$ would be optimal. Indeed, let $\xi$ be any real number such that $\beu < 2$, and let $\eps$ be such that $1 - 1/\beu < \eps < 1/2$. Let $\al $ satisfy $\eps > 1 - 1/\al$, with $\beu < \al < 2$. Then we have $L(\as_i) \leq L(\as_{i-1}) \ll A_i^{-1/\al}$ hence $| \det(\as_{i-1}, \as_i , \as_{i+1}) | \ll A_{i+1} A_i^{-2/\al}$. For $i \in \ids$, this implies $A_{i+1} \gg A_i^{2/\al}$ hence $a_{i+1,0} > 2 a_{i,0}$ for $i$ sufficiently large. Let $\xs_i = \as_i + \as_{i-1}$, for $i \in \ids$. Then $\xs_i $ is primitive (otherwise $\vect_\Q(\as_i, \as_{i-1}) \cap \Z^3$ would strictly contain the subgroup generated by $\as_i$ and $\as_{i-1}$), and $a_{i,0} < x_{i,0} < 2a_{i,0} < a_{i+1,0}$. Therefore $\xs_i$ is not collinear to any minimal point; however it satisfies $L(\xs_i) \leq 2 L(\as_{i-1}) \ll A_i^{-1/\al} \ll X_i^{-1/\al}$ hence $\xs_i \in \ensa_\eps$ if $i \in \ids$ is sufficiently large.

\bigskip

\Dem of Proposition \ref{prop36}:  We may assume $x_0 > 0$. 
Let us denote by $\as_i$ the minimal point corresponding to $x_0$; we have 
$a_{i,0} \leq x_0 < a_{i+1,0}$ and $L(\as_i) \leq L(\xs)$.
Let  $\al > \beu $ be such that $\al + \eps < 2$. 
Using Equation \eqref{eq321} we see that 
$$\gdx L(\as_i) L(\as_{i+1}) \ll \gda_{i+1} L(\as_i) ^2 \ll L(\as_i)^{2-\al}$$
and 
$$\gda_{i+1} L(\as_i) L(\xs) \ll L(\as_i)^{1-\al} \gda_i ^{\eps-1} \ll \gda_i ^{\al + \eps -2}.$$
Therefore Lemma \ref{lem23} proves that the integer $\det(\xs,\as_i,\as_{i+1})$
is zero for $i$ sufficiently large in terms of $\eps$ and $\xi$. This implies $\xs \in \vk$, where $k$ is the index such that $i_k \leq i < i_{k+1}$. 

Let us assume that $\xs$ is not collinear to $\as_i$. Then we have 
$$\hdevk \leq \norme{\xs \wedge \as_i} \ll \gdx L(\xs) 
\ll \gda_{i+1} ^\eps \ll \hdevk^{\eps/(2-\al)}$$
using Equation \eqref{eq24} and Lemma \ref{lem33}. As $\al + \eps < 2$, this is impossible for $k$ 
sufficiently large (in terms of $\eps$ and $\xi$). We obtain in this way $\xs = \lambda \as_i$ for some non-zero integer $\lambda$, which concludes the proof.

\subsection{Linear Independence Properties} \label{subsec41}

The following proposition is very useful as soon as $\beeps < 2$ for some $\eps$ (see 
Equation \eqref{eqscan} above):

\begin{Prop} \label{prop41nv} 
Let $(\us_i)$ be a sequence of minimal points, ordered according to their norms $U_i$,  such that 
$$\limsup_{i \to \infty} \frac{\log U_{i+1}}{- \log L(\us_i)} < 2.$$
Then:
\begin{enumerate}
\item Up to a finite number of terms, $(\ds_k)$ is a subsequence of $(\us_i)$.
\item For any $i$ sufficiently large, $\us_i$ is among the $\ds_k$ if, and only if, $\us_{i-1}$, $\us_i$ and $\us_{i+1}$ are linearly independent.
\end{enumerate}
\end{Prop}

To prove this proposition, we shall use the following lemma:

\begin{Lemme} \label{lem43}
Let $(\us_i)$ be as in Proposition \ref{prop41nv}, and $i$ be   sufficiently large. Denote by $n$ and $m$ the integers such that $\us_i = \as_n$ and $\us_{i+1} = \as_m$. Then the vectors $\as_n$, $\as_{n+1}$, \ldots, $\as_{m-1}$, $\as_m$ belong to a common plane. 
\end{Lemme}

\Dem of Lemma \ref{lem43}: Let  $\beta$ be such that $\limsup \frac{\log U_{i+1}}{- \log L(\us_i)}  < \beta < 2$, and $t$ be such that $n< t < m$. Then we have
$$\vert \det(\as_{t-1}, \as_t, \as_{t+1})\vert \ll A_{t+1} L(\as_{t-1})^2 \leq U_{i+1} L(\us_i)^2
\ll L(\us_i) ^{2 - \beta}$$
so the integer $ \det(\as_{t-1}, \as_t, \as_{t+1})$ is zero if $i$ is sufficiently large. This concludes the proof of Lemma \ref{lem43}.

\bigskip

\Dem of Proposition \ref{prop41nv}: Let $k$ be  sufficiently large, with $\ds_k = \as_p$. If $\ds_k$ were not among the $\us_i$, Lemma \ref{lem43} would apply with $n < p < m$: $\as_{p-1}$, $\as_p$ and $\as_{p+1}$ would be linearly dependent. This is impossible, so the first part of  Proposition \ref{prop41nv} is proved.

Let us prove the second part. If $\us_i = \as_p$ is not among the $\ds_k$, with $i$ sufficiently large, then 
$\as_{p-1}$, $\as_p$ and $\as_{p+1}$  belong to a common plane.  Lemma \ref{lem43} (applied twice) proves that this plane contains also $\us_{i-1}$ and $\us_{i+1}$. To prove the converse statement, assume there is a plane that contains $\us_{i-1}$, $\us_i$  and $\us_{i+1}$, with $\us_i = \as_p$ and  $i$ sufficiently large. Then this plane contains also $\as_{p-1}$ and $\as_{p+1}$, by applying Lemma \ref{lem43} twice; so $\us_i$ is not among the $\ds_k$. This concludes the proof of Proposition \ref{prop41nv}.

\section{The Main Diophantine   Result}  \label{sec5}

This Section is devoted to the statement (\S \ref{subsec51}) and proof (\S \ref{subsec52}) of our main diophantine result. This theorem will be the key point in the proofs of Section \ref{sec6}. 

\subsection{Statement of the Results} \label{subsec51}

Throughout this Section, we define $\eps_1$ by 
\begineq \label{eq512} 
\eps_1 = \Big(2-\beu \Big)  \Big( 2 - \beu  + (2-\bez )\beu \Big). 
\eneq
Since $\beu \leq \bez$, we have  $\eps_1 \geq   (2  - \bez)^2 (1 +\bez)$. 
For instance, if $\bez = 1 + \frac{\sqrt2}{2}$ (denoted by $\sigma_2$ in the Introduction) then $\eps_1 \geq 0.232 $. 

This is the number $\eps_1$ referred to in  Theorem  \ref{thcor516} stated in the Introduction, which is contained in the following result (recall that $\calf$ and $\calr$ were defined in \S \ref{subsec26}, and $\ensa_\eps$ in \S \ref{subsec22}):

\begin{Th} \label{th51} Let $\xi$ be an irrational  non-quadratic real number. 

\begin{enumerate}
\item[(a)] Suppose first $\bez < 2$. Choose a real number $\eps$ such that $0 < \eps < \eps_1$, and let $(\us_i)_{i \geq 1} $ be the sequence of minimal points in $\ensa_{\eps}$, ordered according to their norms $U_i$. Then, we have 
$$L(\us_i) = U_i^{-1+\odu}, \hspace{0.8cm}
\liminf_{i \to \infty} \frac{\log U_{i+1}}{\log U_i} > 1  \hspace{0.6cm} \mbox{and} \hspace{0.6cm}
\bez = \beeps = \limsup_{i \to \infty} \frac{\log U_{i+1}}{\log U_i}.$$
Moreover, there exists a function $\psi \in \calf$ such that $[\us_i, \us_i, \us_{i+1}]$ is collinear to $\us_{\psi(i)}$ for any $i$ sufficiently large. At last, $(\us_i)$ (modulo a finite number of terms) and $\psi$ (modulo $\calr$) are independent from the choice of $\eps$. 
\item[(b)] Conversely, let $(\vs_i)_{i \geq 1}$ be any sequence of primitive points in $\Z^3$ such that:
\begin{itemize}
\item The sequence of first coordinates $(v_{i,0})_{i \geq 1}$ is positive and increasing,
\item As $i \to \infty$, $L(\vs_i)  = V_i^{-1+\odu}$,
\item The real number $\beta = \limsup_{i \to \infty} (\log V_{i+1})/(\log V_i)$ satisfies $\beta < 2$,
\item For any $i$ sufficiently large, if $\vs_{i-1}$, $\vs_i$ and $\vs_{i+1}$ are linearly dependent then $[\vs_i, \vs_i, \vs_{i-1}]$ is collinear to $\vs_{i+1}$.  
\end{itemize}
Then $\bez = \beta < 2$, and $(\vs_i)$ is exactly (up to a finite number of terms) the sequence $(\us_i)$ of Part (a).
\end{enumerate}
\end{Th}

In Part $(a)$, the assertion on  Roy's bracket is equivalent to $\us_{i+1}$ collinear to $[\us_i, \us_i, \us_{\psi(i)}]$ thanks to  \eqref{eq211}. Therefore knowing $\psi$ enables one to construct the sequence $(\us_i)$ by induction.  On the other hand, it follows from Proposition \ref{proplem61} (proved in \S \ref{subsec61} below) that we have
$$\bez = \densi(\psi),$$
where $\densi(\psi)$ was defined in \S \ref{subsec26}. However, knowing $\bez$ does not determine $\psi $ (modulo $\calr$). Indeed, there are functions $\psi$ and $\psi'$, distinct modulo $\calr$, such that $\densi(\psi) = \densi(\psi')$ (see the end of \S \ref{subsec26}). The numbers $\xi$, $\xi'$ constructed from them using the palindromic prefix method  are such that $\bez = \bezpr$ (thanks to Theorem \ref{th63}), but their approximants satisfy completely different recurrence relations. 

\bigskip

In Part $(b)$, the last assumption on $(\vs_i)$ is necessary. Indeed, let us consider the real  number $\xi$ constructed using  the palindromic prefix method (as in  \cite{Alloucheetal} and  \cite{BL}) from the characteristic Sturmian word with slope $[0, 3, 3, 3, \ldots]$. Then $\bez =  1.767\ldots < 2$, but the sequence $(\vs_i)$ consisting in all points denoted by $\ds_k$ and $\es_k$ in \S \ref{subsec52} below satisfies the first three assumptions with $\beta = 1.868\ldots$. 

\bigskip

Theorem \ref{th51}  can be generalized to the case where we assume $\beeps < 2$ for some very small $\eps$. This leads to Th\'eor\`eme 2.2 announced in \cite{CRASasim}; the proof follows the same lines as the one given here. 

\subsection{Proof of Theorem \ref{th51}} \label{subsec52}

The proof falls into eight steps. We let $\xi$ be a real  number such that $[\Q(\xi) : \Q] \geq 3$ and $\bez < 2$. These are the only assumptions throughout the first seven steps. The notation of Theorem \ref{th51} comes into the play only in Step 8. 

Recall that $(\as_i)$ denotes Davenport-Schmidt's sequence of minimal points constructed from $\xi$, and $(\ds_k)$ is the subsequence consisting in all $\as_i$ such that $\as_{i-1}$, $\as_i$ and $\as_{i+1}$ are linearly independent (see \S \ref{subsec23}). As in \S \ref{subsec22}, for $\eps \in \zuouvz$ we let 
$\ensa_\eps$ be the set of all $\xs \in \Z^3 \moins \{  0 \}$ such that $L(\xs) \leq \gdx^{\eps-1}$, where $X =  \max(|x_{0}|, |x_{1}|, |x_{2}|) $.

\smallskip

The sketch of the proof is as follows. We first construct a sequence $(\es_k)$ of minimal points, such that (essentially) $\es_k$ is the first very precise minimal point after $\ds_k$. We have $D_k < E_k \leq D_{k+1}$, and we show how $\ds_k$, $\es_k$ and $\hdevk$ are interrelated. Then we construct a sequence $(\ws_t)$ of minimal points, which turns out to be the sequence of ``all'' very precise approximants; $(\ds_k)$ and $(\es_k)$ are subsequences of this sequence, and other points between $\es_k$ and $\ds_{k+1}$ are constructed thanks to an interpolation procedure using Roy's bracket. The key point of the proof is the definition and properties of this sequence $(\ws_t)$, and of the function $\psi$ associated with it. Then Theorem \ref{th51} follows easily : the sequences $(\us_i)$ and $(\vs_i)$ in Theorem \ref{th51} are proved to coincide (up to a finite number of terms) with $(\ws_t)$.

\medskip

\begin{Step} \label{step1} Construction of  the sequence $(\es_k)$.
\end{Step}

Let us fix a real number $\eps_2$ such that $0 < \eps_2 < \eps_1$, where $\eps_1$ is defined by Equation \eqref{eq512}. Let $(\bs_t)_{t \geq 1}$ be the sequence of all minimal points in $\ensa_{\eps_2}$, ordered according to their norms $B_t$. Since $\beepsde \leq \bez < 2$, 
Proposition \ref{prop41nv} applies (thanks to \eqref{eqscan}) and shows that for any $k$ sufficiently large $\ds_k$ is equal to some $\bs_t$; then we let $\es_k = \bs_{t+1}$. 

The sequence $(\es_k)_{k \geq 1}$ defined in this way (by choosing arbitrarily $\es_k$ for small values of $k$) seems to depend on the choice of $\eps_2$, but actually this dependence concerns only finitely many terms,  as  the following lemma shows.

\begin{Lemme} \label{lem118}
Let $\eps$ be such that $0 < \eps \leq \eps_2$. Then for any $k$ sufficiently large (in terms of $\eps$), $\ds_k$ and $\es_k$ are consecutive elements of the sequence of minimal points in $\ensa_{\eps}$. 
\end{Lemme}

It is possible to state more precise versions of this lemma (see Lemma \ref{lem120} below). We shall prove the following one now:

\begin{Lemme} \label{lem118bis}
Let  $\beta_V \in (1,2)$, and   $(\vs_i)_{i \geq 1}$ be any sequence of minimal points in $\ensa_{\eps_2}$, ordered according to their norms $V_i$,  such that 
$$\limsup_{i \to \infty} \frac{\log V_{i+1}}{- \log L(\vs_i)} \leq \beta_V.$$
Assume that for some $\eps' > 0$ with  
$$\eps' < (2 - \beu)(2 - \beu + (2 - \beta_V)\beu), $$
we have $\es_k \in \ensa_{\eps'}$ when $k$ is sufficiently large. 

Then, for any $k$ sufficiently large, $\ds_k$ and $\es_k$ are consecutive elements of the sequence $(\vs_i)_{i \geq 1}$. 
\end{Lemme}

To deduce Lemma \ref{lem118}, it suffices to apply Lemma \ref{lem118bis} to the sequence $(\vs_i)$ of all minimal points in $\ensa_{\eps}$, with $\eps' = \eps_2 < \eps_1$ and $\beta_V = \beeps \leq \bez$ thanks to Equation \eqref{eqscan}. The full generality of Lemma \ref{lem118bis} will be useful in Step \ref{step8}, since it will turn out that $\es_k \in \ensa_{\eps'}$ for any $\eps' > 0$ (as soon as $k$ is sufficiently large in terms of $\eps'$). 

\bigskip

\Dem of Lemma \ref{lem118bis}: 
Let $k$ be sufficiently large. Proposition \ref{prop41nv} provides an integer $i$ such that $\ds_k = \vs_i$. Since  $(\vs_i)$ is a subsequence of $(\bs_t)$,  we have $V_{i+1} \geq E_k$. Let us assume that $V_{i+1} > E_k$. Since $\ds_{k+1}$ is among the $\vs_i$, we have $D_k < E_k < V_{i+1} \leq D_{k+1}$. Choosing $\beta \in (\beta_V, 2)$, we have $V_{i+1} \ll L(\vs_i)^{-\beta}$ hence, by Lemma \ref{lem3131}:
$$\hdevk \ll E_k L(\ds_k) \mbox{ and } \hdevk \ll V_{i+1} L(\es_k) \ll L(\ds_k)^{-\beta} L(\es_k).$$
The  product of these relations yields, by choosing $\al \in (\beu, 2)$ and using Lemmas \ref{lem33} and \ref{lem34}:
\begin{eqnarray*}
\hdevk ^2 
&\ll& L(\ds_k)^{1-\beta} E_k^{\eps'}\\
&\ll& D_k ^{\beta-1} D_{k+1} ^{\eps'} \ll \hdevk^{\al(\beta-1)+\frac{\eps'}{2-\al}}
\end{eqnarray*}
hence $\eps' \geq (2-\al)(2 - \al(\beta-1))$. This contradicts the assumption on $\eps'$ if $\al$ and $\beta$ are close enough to $\beu$ and $\beta_V$.

\medskip

\begin{Step} \label{step2} Relations between $\ds_k$ and $\es_k$.
\end{Step}

The following properties will be used many times in the proof of Theorem \ref{th51}, sometimes without reference:
\begin{eqnarray}
&L(\ds_k) = D_k^{-1+\odu}, \hspace{0.6cm} L(\es_k) = E_k^{-1+\odu} \label{eq41}\\
&D_k < E_k \leq D_{k+1}, \hspace{0.6cm} E_k \leq D_k ^\beta \label{eq42}\\
&\hdevk D_k^{1+ \odu}  \leq E_k \label{eq43}\\
&D_k^{1/\al} \leq \hdevk , \hspace{0.6cm} D_{k+1}^{2-\al}  \leq \hdevk .  \label{eq44}
\end{eqnarray}
In these formulas, $\al$ and $\beta$ are chosen such that $\beu < \al < 2$ and $\bez < \beta < 2$; in the applications, they  are assumed to be sufficiently close to $\beu$ and $\bez$. The inequalities hold for $k$ large enough (in terms of $\al$ and $\beta$). The symbol $\odu$ denotes a sequence (possibly depending on $\al$ and $\beta$) that tends to zero as $k$ tends to infinity. 

Equations \eqref{eq41} and \eqref{eq42} follow immediately from Lemma \ref{lem118}. Thanks to Lemma \ref{lem3131}, they imply \eqref{eq43}. At last, Equation \eqref{eq44} follows from Lemmas \ref{lem33} and \ref{lem34}. 

We will use repeatedly the following consequence of these relations:
$$D_k^{\odu}= E_k^{\odu}= D_{k+1}^{\odu}= \hdevk^{\odu}.$$

\medskip

\begin{Step} \label{step3} Construction of a sequence $(\ws_t)$ of primitive points.
\end{Step}

Let $k$ be a sufficiently large integer, say $k \geq k_0$. 

We define now an integer $s_k \geq 1$, and points $\ns_0$, \ldots, $\ns_{s_k} \in \Z^3$, as follows. To begin with, let $\ns_0 = \ds_{k+1}$. For any $\sigma \geq 0$ such that $\ns_\sigma$ is defined and $\gdn_\sigma \geq E_k$, we let $\ns_{\sigma+1}$ be the primitive element of $\Z^3$, with non-negative first coordinate, collinear to $[\ns_\sigma, \ds_k, \es_k]$. At last, we let $s_k$ be the greatest integer $\sigma$ for which $\ns_\sigma$ is defined (and $s_k = + \infty$ if $\ns_\sigma$ is defined for any $\sigma$). 

We let $(\ws_t)_{t \geq 1}$ be the sequence of all points $\ns_\sigma$, with $k \geq k_0$ and $0 \leq \sigma \leq s_k$, ordered according to their norms $W_t$. Since $\det(\ds_k)$ and $\det(\es_k)$ are non-zero for $k$ sufficiently large, we have $\det (\ws_t) \neq 0$ for $t$ sufficiently large.

\medskip

\begin{Step} \label{step4} First properties of the  sequence $(\ws_t)$.
\end{Step}

In this step, we prove the following lemma:

\begin{Lemme} \label{lem119}
We have 
\begineq \label{eq522}
s_k \leq \frac{1 + \odu}{2 - \beu}  \mbox{ as } k \to \infty, 
\eneq
\begineq \label{eq18bis}
L(\ns_\sigma) = {\gdn_\sigma }^{-1+\odu}, 
\eneq
and
\begineq \label{eq521}
\left\{
\begin{array}{l}
\mbox{all $\ns_\sigma$ are minimal points, with $D_{k+1} = \gdn_0 > \ldots > \gdn_{s_k} = D_k$;}\\
\mbox{in particular $\ns_{s_k} = \ds_k$ and $\ns_{s_k-1} = \es_k$. }
\end{array}
\right.
\eneq
Accordingly, $s_k$ is finite. 
\end{Lemme}

\Demdeuxpoints Thanks to Equation \eqref{eq41}, Proposition \ref{prop21} yields, for $k \geq k_0$ and $\sigma \in \unsk$: 
\begineq \label{eq526}
\gdn_\sigma \leq \gdn_{\sigma-1} D_k E_k^{-1+\odu} \mbox{ and }
L(\ns_\sigma) \leq L(\ns_{\sigma-1}) D_k ^{-1+\odu} E_k.
\eneq
By induction on $\sigma$, this implies 
\begineq \label{eq523}
\gdn_\sigma \leq D_{k+1} \Big( D_k E_k ^{-1+\odu} \Big) ^\sigma
\eneq
and
\begineq \label{eq524}
L(\ns_\sigma) \leq D_{k+1}^{-1+\odu} \Big( D_k^{-1+\odu} E_k \Big) ^\sigma
\eneq
for $\sigma \in \zerosk$, 
where the sequences involved in the notation $\odu$ tend to 0 as $k$ tends to infinity, uniformly with respect to $\sigma$. 
Since $L(\ns_\sigma) \gdn_\sigma \gg 1$, it turns out that \eqref{eq523} and \eqref{eq524} are actually equalities.

Moreover, \eqref{eq43}, \eqref{eq44} and \eqref{eq526} imply 
\begineq \label{eq529}
\gdn_\sigma < \gdn_{\sigma-1}
\eneq
if $k$ is large enough and $1 \leq \sigma \leq s_k$. 

\medskip

Let us prove \eqref{eq522} now. Let $k \geq k_0$ and $\sigma \in \zeroskmu$. Using   Equations \eqref{eq44} and  \eqref{eq523},  and the fact that $\gdn_\sigma \geq E_k$, we obtain:
$$\hdevk^{1/(2-\al)} \gg D_{k+1} \geq D_k ^{-\sigma} E_k ^{(\sigma+1)(1+\odu)}.$$
Now Equation \eqref{eq43} yields 
\begineq \label{eq529bis}
\hdevk^{1/(2-\al)} \geq \hdevk ^{(\sigma+1)(1+\odu)} D_k^{1 + \sigma \odu}.
\eneq
If \eqref{eq522} fails to hold then for infinitely many such $\sigma$ then we have $(\sigma+1)(1 + \odu) > 1/(2-\al)$ (if $\al$ is chosen small enough), hence \eqref{eq529bis} yields $1 + \sigma \odu < 0$ and (using Equation \eqref{eq44})
$$\frac{1}{2 - \al} \geq (\sigma + 1)(1 + \odu) + \al (1 + \sigma \odu), $$
which implies
$$(\sigma+1)(1 + \odu) \leq \frac{1}{2-\al} - \al.$$
This contradiction concludes the proof of \eqref{eq522}. In particular, $\sigma$ is bounded uniformly with respect to $k$, hence $\sigma \odu = \odu$. Therefore the following inequality follows from \eqref{eq44}, \eqref{eq523},  and \eqref{eq524}:
\begineq \label{eq525} 
\gdn_\sigma L(\ns_\sigma) \leq \hdevk^{\odu} \mbox{ for any } \sigma \in \zerosk. 
\eneq
Now $\hdevk \leq E_k D_k^{-1+\odu}  \leq E_k \leq \gdn_\sigma$ if $\sigma \in \zeroskmu$. In this case, Equation \eqref{eq525} and Proposition \ref{prop36} imply that $\ns_\sigma$ is a minimal point. Let us conclude the proof of \eqref{eq521} now.

We know that $\ns_{s_k-1}$ is a minimal point, with $\gdn_{s_k-1} \geq E_k$ hence $L(\ns_{s_k-1}) \leq L(\es_k)$. Equations \eqref{eq41}, \eqref{eq42}  and \eqref{eq526} yield
\begineq \label{eq528} 
L(\ns_{s_k}) \leq D_k^{-1+\odu}. 
\eneq
 If $\gdn_{s_k} < D_k$ then \eqref{eq528} gives $\ns_{s_k} \in \ensa_{\odu}$; otherwise we have $\gdn_{s_k} \geq D_k \gg E_k^{1/2} \geq \hdevk^{1/2}$ and \eqref{eq525} implies 
$\ns_{s_k} \in \ensa_{\odu}$. Therefore Equation \eqref{eq18bis} holds in both cases. Then  Proposition \ref{prop36} proves that $\ns_{s_k}$ is a minimal point, say $\as_i$. If $i \leq i_k-1$ (that is, $\gdn_{s_k} < D_k$) then Equations \eqref{eq528} and \eqref{eq44}  yield
$$D_k ^{2-\al} \ll \hdevkmu \equ D_k L(\as_{i_k-1}) \leq D_k L(\ns_{s_k}) \leq D_k^{\odu}$$
which is impossible. Therefore $\gdn_{s_k} \geq D_k$; but we have also $\gdn_{s_k} < E_k$ by definition of $s_k$. So $\ns_{s_k}$ is a minimal point in $\ensa_{\odu}$ between $\ds_k$ and $\es_k$, distinct from $\es_k$. Lemma \ref{lem118}  proves that $\ns_{s_k} = \ds_k$. 

By definition of $\ns_{s_k}$, this proves that $\ds_k$ is  collinear to $[\ns_{s_k-1}, \ds_k, \es_k]$. Equation \eqref{eq211} implies that $\es_k$ is collinear to $[\ns_{s_k-1}, \ds_k, \ds_k]$, hence to $\ns_{s_k-1}$ thanks to \eqref{eq212}. This concludes the proof of Lemma \ref{lem119}.

\medskip

\begin{Step} \label{step5} Connection between $(\ws_t)$ and $\ensa_\eps$.
\end{Step}

The following result is a stronger version of Lemma \ref{lem118}:

\begin{Lemme} \label{lem120}
Let $\eps \in (0,  \eps_1)$, and $(\us_i)$ be the sequence of all minimal points in $\ensa_\eps$, ordered according to their norms $U_i$. Then the sequences $(\us_i)$ and $(\ws_t)$ coincide up to a finite number of terms.
\end{Lemme}

\Demdeuxpoints Equation \eqref{eq18bis} means $L(\ws_t) = W_t^{-1+\odu}$, and shows that $\ws_t$ belongs to $\ensa_{\eps}$ for $t$ sufficiently large in terms of $\eps$. Assume there is a minimal point $\as_j$ in $\ensa_{\eps}$, with $j$ arbitrarily large,  such that $W_t < A_j < W_{t+1}$. Then we can write $\ws_t = \ns_{\sigma+1}$ and $\ws_{t+1} = \ns_{\sigma}$ for some $k $ and some $\sigma \in \{0, \ldots, s_k-1\}$. Lemma \ref{lem3131} yields
$$\hdevk \ll A_j L(\ns_{\sigma+1}) \mbox{ and }
\hdevk \ll \gdn_{\sigma} L(\as_j)$$
hence, by choosing $\al \in (\beu, 2)$ and $\beta \in (\bez, 2)$ and  using \eqref{eq523}, \eqref{eq524}, \eqref{eq42}, and \eqref{eq44}:
$$\hdevk^2 \ll A_j^\eps E_k D_k^{-1+\odu} \ll D_{k+1}^\eps D_k^{\beta-1+\odu} \ll
\hdevk^{\frac{\eps}{2-\al}+\al(\beta-1) + \odu}.$$
Since $\eps < \eps_1$, this gives a contradiction when $\al$ and $\beta$ and sufficiently small, and $k$ is sufficiently large. This concludes the proof of Lemma \ref{lem120}.

\medskip

\begin{Cor} \label{cor121}
We have
$$L(\ws_t) = W_t^{-1+\odu}$$
and
$$\bez = \beeps = \limsup_{t \to \infty} \frac{\log W_{t+1}}{\log W_t}$$
for any $\eps \in (0, \eps_1)$. 
\end{Cor}

\Demdeuxpoints The first equality is nothing but Equation \eqref{eq18bis}. The second one follows from it, thanks to Equation \eqref{eqscan} and Lemma \ref{lem120}.

\medskip

\begin{Step} \label{step6} Linear dependence between consecutive terms of the sequence $(\ws_t)$.
\end{Step}

\begin{Lemme} \label{lem122}
Let $t$ be a sufficiently large integer such that $\ws_{t-1}$, $\ws_t$ and $\ws_{t+1}$ are linearly dependent. Then the bracket $[\ws_t, \ws_t, \ws_{t+1}]$ is collinear to $\ws_{t-1}$, so that $[\ws_t, \ws_t, \ws_{t-1}]$ is collinear to $\ws_{t+1}$.
\end{Lemme}

\Demdeuxpoints Thanks to Corollary \ref{cor121}, Proposition \ref{prop41nv} applies (since $\bez  < 2$) and provides an integer $k$ such that $D_k \leq W_{t-1} < W_t < W_{t+1} \leq D_{k+1}$. By construction, there exist non-zero integers $\lambda$, $\lambda'$ such that 
$$\lambda \ws_t = [\ws_{t+1}, \ds_k, \es_k] = - \ws_{t+1} J \es_k J \ds_k $$
and 
$$\lambda ' \ws_{t-1} = [\ws_{t}, \ds_k, \es_k] = - \ws_{t} J \es_k J \ds_k .$$
This implies 
\begin{eqnarray*} 
[\ws_t, \ws_t, \ws_{t+1}]
&=& -\ws_t J \ws_{t+1} J \ws_t \\
&=&   \lambda^{-1} \ws_t (J \ws_{t+1})^2 J \es_k J \ds_k \\
&=& \det(\ws_{t+1}) \lambda^{-1} \lambda' \ws_{t-1}
\end{eqnarray*}
since $(J \xs)^2 = - \det(\xs) \Id$ for any $\xs \in \Z^3$ identified with a symmetric matrix.
This concludes the proof of the first statement of Lemma \ref{lem122}; the second one immediately follows from Equation \eqref{eq211}.

\medskip

\begin{Step}  \label{step7} Construction and properties of $\psi$.
\end{Step}

Let $\beta $ be such that $\bez < \beta < 2$, and $\capa = \beta/2$. Corollary  \ref{cor22} proves, thanks to Corollary \ref{cor121}, that $[\ws_t, \ws_t, \ws_{t+1}] \in \ensa_{\odu}$. By Proposition \ref{prop36}, this bracket is collinear to a minimal point for $t$ sufficiently large; and by Lemma \ref{lem120} this minimal point belongs to the sequence $(\ws_t)$. 

Therefore for any sufficiently large integer $t$, there exists a unique integer, denoted by $\psi(t)$, such that $[\ws_t, \ws_t, \ws_{t+1}] $ is collinear to $\ws_{\psi(t)}$. If $t$ is not sufficiently large, we let $\psi(t) = t-1$.

In this definition, $\psi$ depends on the choice of a threshold  for determining when $t$ is ``sufficiently large'',   on the integer $k_0$ chosen at the beginning of Step \ref{step3} (in the sequence $(\ws_t)$, the index $t$ is shifted if the choice of $k_0$ is modified), and on finitely many arbitrary choices of initial values (see Step \ref{step1}). 
Actually this dependence is very mild, as shows the following lemma.

\begin{Lemme} \label{lem123}
\begin{itemize}
\item The function $\psi$ belongs to the set $\calf$ defined in \S \ref{subsec26}, and up to $\calr$ it   depends only upon $\xi$.
\item For any $t$ sufficiently large, the following assertions are equivalent:
	\begin{itemize}
	\item[$(i)$] $\psi(t) = t-1$.
	\item[$(ii)$] $\ws_{t-1}$, $\ws_t$ and $\ws_{t+1}$ are linearly dependent.
	\item[$(iii)$] $\ws_t$ is not among the $\ds_k$.
	\end{itemize} 
\end{itemize}
\end{Lemme}

\Demdeuxpoints 
Let us start with the second statement. Thanks to Corollary \ref{cor121},  Proposition \ref{prop41nv} applies and proves that $(ii) \Leftrightarrow (iii)$. Since  $\ws_{t}$, $\ws_{t+1}$ and  $[\ws_t, \ws_t, \ws_{t+1}]$ are always linearly dependent, the implication $(i) \Rightarrow (ii)$ is obvious. At last, the implication  $(ii) \Rightarrow (i)$ follows from Lemma \ref{lem122}. 

\smallskip

Let us prove the first statement now. Let $t$ and $k$ be sufficiently large integers such that $\ws_t = \ds_k$. Proposition \ref{prop21}, Corollary \ref{cor121}, and Equations \eqref{eq43} and \eqref{eq44} yield, since $\ws_{t+1} = \es_k$:
$$W_{\psi(t)} \ll W_t^2 W_{t+1}^{-1+\odu} \ll W_{t} ^{1+\odu} \hdevk^{-1} \ll W_t^{1 - 1/\al + \odu}$$
hence $\psi(t) \leq t-1$. As $\psi$ is defined by $\psi(t) = t-1$ for small values of $t$, this proves that $\psi(t) \leq t-1$ for any $t$.

Let $(\teta_k)_{k \geq 0}$ be the sequence associated with $\psi$, as in \S \ref{subsec26}; namely, $(\teta_k)$ is the increasing sequence of all indices $t$ such that $\psi(t) \leq t-2$. Since $(\ds_k)$ is a subsequence of $(\ws_ t)$, the second part of Lemma \ref{lem123} shows that 
$$\ws_{\teta_k} = \ds_k \mbox{ for any } k$$ 
up to shifting the index $k$; in particular, the sequence $(\teta_k)$ is infinite. 

Let us prove that $\psi$ is asymptotically reduced (as defined in  \S \ref{subsec26}, that is $\psi(\teta_k) <  \teta_{k-1}$ and $\psi(\teta_k) \neq  \psi(\teta_{k-1})$ for any $k$ sufficiently large). 

First, assume that $\psi(\teta_k) \geq \teta_{k-1}$ with $k$ sufficiently large, hence $W_{\psi(t)} \geq D_{k-1}$ with $t = \teta_k$ (i.e., $\ws_t = \ds_k$). Then $\ws_{\psi(t)}$ is a minimal point between $\ds_{k-1}$ and $\ds_k = \ws_t$, so that $\ws_{\psi(t)} \in \vkmu = \Span_{\Q} (\ds_{k-1}, \ds_k) \cap \Z^3$. But $\ws_{\psi(t)}$ is collinear to $[\ws_t, \ws_t, \ws_{t+1}]$, so that $\ws_{\psi(t)} \in \vk = \Span_{\Q} (\ws_{t}, \ws_{t+1}) \cap \Z^3$. Now $\vk \cap \vkmu = \Z \ds_k$ by definition of the sequence $(\ds_k)$, therefore $\ws_{\psi(t)}$ is collinear to $\ds_k = \ws_t$: this is impossible. So we have  $\psi(\teta_k) < \teta_{k-1}$ for any $k$  sufficiently large.

Let us assume now that $\psi(\teta_k) =  \psi(\teta_{k-1})$ with $k$ sufficiently large. Then $\ws_{\psi(\teta_k)} = \ws_{\psi(\teta_{k-1})}$ is proportional to both $[\ds_k, \ds_k, \es_k]$ and $[\ds_{k-1}, \ds_{k-1}, \es_{k-1}]$ so it belongs to $\vkmu \cap \vk = \Z \ds_k$. This is impossible since $\psi(\teta_k) \leq \teta_k - 1$, so that  $\psi(\teta_k) \neq   \psi(\teta_{k-1})$ for $k$ sufficiently large. 

This concludes the proof that $\psi$ is asymptotically reduced. Using Corollary \ref{cor121} (see also Lemma A.2 of \cite{SFcombi}), it is not difficult to deduce from the property $\bez < 2$ the existence of a $c$ such that $\psi(t) \geq t-c$ for any $t$. So we have proved that $\psi \in \calf$.

As noticed before the statement of  Lemma \ref{lem123}, $\psi$ depends only on the choice of some initial values. If these choices were different, then  $\psi(n)$ would be replaced (for $n$ sufficiently large) with $\psi(n+c) + c'$ for some integers $c$, $c'$ independent from $n$.   This means exactly that $\psi$ remains the same modulo $\calr$, thereby concluding  the proof of Lemma \ref{lem123}.

\medskip

\begin{Step}  \label{step8} Proof of Theorem \ref{th51}.
\end{Step}

Let us prove Part $(a)$ first. Let $\eps \in (0, \eps_1)$, and $(\us_i)$ be the sequence of all minimal points in $\ensa_\eps$, ordered according to their norms $U_i$. Up to shifting indices, Lemma \ref{lem120} shows that $\us_i = \ws_i$ for any $i$ sufficiently large (and the value of $\eps$ has an influence only on the shift). So Corollary \ref{cor121} implies  $L(\us_i) = U_i^{-1+\odu}$ and 
$\bez = \beeps = \limsup_{i \to \infty} \frac{\log U_{i+1}}{\log U_i}$. 

Let $i$ be sufficiently large, and $k$ be such that $D_k \leq W_{i-1} < D_{k+1}$. Since \eqref{eq523} is an equality, we have (using \eqref{eq43} and \eqref{eq44}):
$$ W_i W_{i-1}^{-1} = E_k^{1 + \odu} D_k ^{-1} \geq \hdevk^{1 + \odu} \geq W_{i-1}^{2 - \al + \odu}$$
for any (fixed) $\al \in (\beu, 2)$, hence
$$\liminf_{i \to \infty} \frac{\log W_i}{\log W_{i-1}} \geq 3 - \beu > 1.$$
Since the function $\psi \in \calf$ has been constructed in Step \ref{step7}, 
this concludes the proof of Part $(a)$. 

\smallskip

Let us prove Part $(b)$ now. Let $(\vs_i)$ and $\beta < 2$ be as in Theorem \ref{th51}. Since $L(\vs_i) = V_i ^{-1+\odu}$, for any $i$ sufficiently large we have $\vs_i \in \ensa_{\eps_2}$, and $\vs_i$ is a minimal point thanks to Proposition \ref{prop36}. 

Let $k$ be sufficiently large. Proposition \ref{prop41nv} shows that $\ds_k = \vs_i$  for some $i$, and Lemma \ref{lem118bis} (applied with a small $\eps ' > 0$ thanks to Corollary \ref{cor121}, and $\beta_V  = \beta$) yields $\es_k = \vs_{i+1}$. Let $t$ be such that $\ws_t = \ds_k$ (that is, $t = \teta_k$ with the notation of Step \ref{step7}). We have $\ws_t = \vs_i$ and $\ws_{t+1} = \vs_{i+1}$; let us prove by induction that $\ws_{t + \ell} = \vs_{i + \ell}$ for any $\ell \in \{ 0, \ldots, \ell_0\}$,  where $\ell_0$ is such that $W_{t + \ell_0} = D_{k+1}$. 

We may assume $\ell_0 \geq 2$, otherwise this is proved already. Let $\ell  \in \{ 1, \ldots, \ell_0 - 1\}$ be such that  $\ws_{t + \ell} = \vs_{i + \ell}$. Since $\ell_0 \geq 2$, we have $D_k \leq W_{t+\ell-1} < W_{t+\ell} < W_{t+\ell+1} \leq D_{k+1}$ so $\ws_{t+\ell}$ is not among the $\ds_k$. Lemma \ref{lem123} yields $\psi(t+\ell) = t+\ell-1$, so that $\ws_{t+\ell+1}$ is the primitive point, with non-negative first coordinate, 
collinear to $[\ws_{t+\ell}, \ws_{t+\ell} , \ws_{t+\ell-1} ]$. By induction hypothesis, this bracket is equal to 
$[\vs_{i+\ell}, \vs_{i+\ell} , \vs_{i+\ell-1} ]$. Now $V_{i+\ell} = W_{t+\ell} < W_{t+\ell+1} \leq D_{k+1}$ yields $V_{i+\ell+1} \leq D_{k+1}$, since $(\ds_k)$ is a subsequence of $(\vs_i)$ up to a finite number of terms. So $\vs_{i+\ell-1}$, $\vs_{i+\ell}$ and $\vs_{i+\ell+1}$ are linearly dependent. By assumption on the sequence $(\vs_i)$, this proves that $\vs_{i+\ell+1}$ is collinear to $[\vs_{i+\ell}, \vs_{i+\ell} , \vs_{i+\ell-1} ]$, and to $\ws_{t+\ell+1}$. Since both vectors are primitive and have  non-negative first coordinates, we have 
$\ws_{t + \ell+1} = \vs_{i + \ell+1}$, thereby concluding the induction.

Since this result holds for any $k$ sufficiently large, we have proved that   the sequence $(\vs_i)$ is obtained from $(\ws_t)$ by shifting the index (up to a finite number of terms). In particular, we have $\beta = \bez$ thanks to Corollary \ref{cor121}, thereby concluding the proof of Theorem \ref{th51}.

\section{The Palindromic Prefix Method} \label{sec6}

Roy's construction \cite{RoyCRAS} of a number $\xi$ such that $\beu = \nbor$ makes a crucial use of a specific word $w$ with many palindromic prefixes, namely the Fibonacci word (see Example \ref{explefibo} above). Bugeaud and Laurent have generalized his construction \cite{BL}  to  any characteristic Sturmian word $w$ (see also \cite{Alloucheetal}). We explain here how to do  it, more generally,  with  any word $w$ such that $\densi(w)  < 2$. 

This method produces a real number $\xi$; the key property (proved in \S \ref{subsec62} below) is that 
\begin{equation} \label{eqmotippm}
\beu \leq \bez \leq \densi(w).
\end{equation}
The main difficulty in proving \eqref{eqmotippm}  is to study the asymptotic behaviour of a sequence defined by induction using Roy's bracket. This was done by Roy and Bugeaud-Laurent by using specific properties of the words $w$ they were considering; we do it in the general case in \S \ref{subsec61}. The proof is elementary, but rather technical.

The difficult diophantine question is to know whether \eqref{eqmotippm} is an equality. Using the results of Section \ref{sec5}, we prove in \S  \ref{subsec62} that 
$$\bez = \densi(w)$$
but we do not know whether the equality $\beu = \bez$ holds for a general word $w$ (it does if $w$ is characteristic Sturmian, as proved by Bugeaud-Laurent). 

\subsection{Asymptotic Behaviour of Some Sequences} \label{subsec61}

The following Proposition (in which $\densi(\psi)$ is defined in \S \ref{subsec26}) 
will be useful in \S \S \ref{subsec62} and \ref{subsec53} below. It is a more general version of the arguments used in \cite{RoyCRAS} (end of the proof of Th\'eor\`eme 2.2), \cite{Roytwoexp} (Lemma 5.2) and \cite{BL} (Lemma 4.1). The proof essentially involves Lemma \ref{lem62} below, which is a generalization of Proposition 6.2 of \cite{SFcombi}. 

\begin{Prop} \label{proplem61}
Let $\xi \in \R$ and $\psi \in \calf$. Let $(\vs_n)_{n\geq 1} $ be a sequence of non-zero points in $\Z^3$ such that $[\vs_n, \vs_n, \vs_{n+1}]$ is collinear to $\vs_{\psi(n)}$ for any sufficiently large $n$, 
\begineq \label{eq611} 
L(\vs_n) = V_n ^{-1 + \odu}
\eneq
and $V_{n+1} \geq V_n ^{\delta}$ for some $\delta > 1$ and any sufficiently large $n$. Then we have 
\begineq \label{eq612} 
\limsup_{n \to + \infty} \frac{\log V_{n+1}}{ \log V_n} = \densi(\psi).
\eneq
\end{Prop}

\Dem of Proposition \ref{proplem61}: For any $n \geq 1$, let $\lambda_n \in \Qetoile$ be such that $[\vs_n, \vs_n, \vs_{n+1}] = \lambda_n \vs_{\psi(n)}$. Proposition \ref{prop21} and Equation \eqref{eq611} yield $|\lambda_n | V_{\psi(n)} \ll V_n^2 V_{n+1}^{-1+\odu}$ and 
$|\lambda_n | L(\vs_{\psi(n)}) \ll V_n^{-2+\odu} V_{n+1}$ hence
$|\lambda_n |^2 V_{\psi(n)} L(\vs_{\psi(n)}) \ll V_n^{\odu} V_{n+1}^{\odu}$.
Therefore \eqref{eq231} yields $|\lambda_n |^2 \ll V_{n+1}^{\odu} $. Now let $a_n, b_n \in \Z$ be coprime integers such that $\lambda_n = a_n / b_n$. Then $\vs_{\psi(n)} \in b_n \Z^3$, so Equations \eqref{eq611} and \eqref{eq231} imply $|b_n| = V_{n+1}^{\odu}$, hence $|a_n| = V_{n+1}^{\odu}$. Finally $|\lambda_n| = V_{n+1}^{\odu}$, and the above inequalities (together with Equation \eqref{eq231}) imply $V_{n+1} \ll V_n^{2+\odu}$ and $V_{\psi(n)} = V_n^{2+\odu} V_{n+1}^{-1}$. Therefore   applying the following lemma   with $u_n = \log V_n$ is enough to  conclude the proof of Proposition \ref{proplem61}. 

\medskip

\begin{Lemme} \label{lem62}
Let $\psi \in \calf$, and $(u_n)_{n \geq 1}$ be a sequence of positive real numbers such that 
$$u_{n+1} = 2 u_n - u_{\psi(n)} + \oduu$$
and $\liminf u_{n+1} / u_n > 1$. Then we have 
\begineq \label{eq623} 
\limsup_{n \to + \infty} \frac{u_{n+1}}{u_n} = \densi(\psi).
\eneq
\end{Lemme}

In the proof of Lemma \ref{lem62}, we shall use the following result:

\begin{Lemme} \label{lem62bis}
Let $\psi \in \calf$, and $(\uti_n)_{n \geq 1}$ be an increasing sequence of positive integers such that 
$$\uti_{n+1} = 2 \uti_n - \uti_{\psi(n)} \mbox{ for any $n$ sufficiently large.}$$
Then we have $\liminf \uti_{n+1} / \uti_n > 1$ and $\limsup \uti_{n+1} / \uti_n  = \densi(\psi)$.
\end{Lemme}

Thanks to this lemma, we see that in Lemma \ref{lem62} the assumption $\liminf u_{n+1} / u_n > 1$ is automatically fulfilled if $u_{n+1} = 2 u_n - u_{\psi(n)} $ for any  $n$ sufficiently large. However, as it stands, Lemma  \ref{lem62} would be false without this assumption since $u_n$ could be a polynomial in $n$. 

\bigskip

\Dem of Lemma \ref{lem62bis}: 
The assertion on the upper limit is exactly Proposition 6.2 of \cite{SFcombi}, already recalled in \S \ref{subsec26}. Let us prove the statement on the lower limit.

Let $(\teta_k)$ be the sequence associated with $\psi$, as in \S \ref{subsec26} (namely, this is the increasing sequence of all integers $n$ such that $\psi(n) \leq n-2$).  For $k$ sufficiently large, we have $\psi(\teta_{k+1}) \leq \teta_k - 1$ (since $\psi \in \calf$) and $\psi(n) = n-1$ for $\teta_k < n < \teta_{k+1}$ hence:
$$\uti_{\teta_{k+1}-1} \geq 2 \uti_{\teta_{k+1}}  - \uti_{\teta_{k}-1} =
\uti_{\teta_{k+1}} + \uti_{\teta_{k+1}-1}  + \uti_{\teta_{k}+1} - \uti_{\teta_{k}} - \uti_{\teta_{k}-1}    $$
(as in the proof of Lemma A.1 of  \cite{SFcombi}).This implies (by induction) the existence of some integer $c$ such that $  \uti_{\teta_{k}+1} -  \uti_{\teta_{k}}  -  \uti_{\teta_{k}-1} \geq c$, hence  
$  \uti_{\teta_{k}+1} \geq \frac54   \uti_{\teta_{k}}$, for any $k$ sufficiently large (since $\uti_n \leq 2 \uti_{n-1}$ for $n$ sufficiently large, and $\uti_n \to \infty$). Now for $\teta_k  < n \leq \teta_{k+1}$ we have
$$\uti_{n} - \uti_{n-1} =  \uti_{\teta_{k}+1} -   \uti_{\teta_{k}} \geq \frac14  \uti_{\teta_{k}} 
\geq  \frac{1}{2^B} \uti_{\teta_{k+1}-1} \geq  \frac{1}{2^B} \uti_{n-1}$$
where $B$ is such that    $\teta_{k+1} - B \leq \psi(\teta_{k+1}) < \teta_{k}$ hence $\teta_{k+1} \leq \teta_{k} + B-1$   for any $k$; such a $B$ exists since $\psi \in \calf$. Therefore we have 
$\uti_n \geq (1  + 2^{-B}) \uti_{n-1}$
for any $n$ sufficiently large, thereby concluding the proof of Lemma \ref{lem62bis}.

\bigskip 

\Dem of Lemma \ref{lem62}: Let $n_0$ be a sufficiently large integer. We let $\eps_n = u_{n+1} - u_n$; then $\eps_n > 0$ for $n \geq n_0$. Let $\delta_n \in \R$ be defined by $u_{n+1} = 2u_n - u_{\psi(n)} + \delta_n u_{n+1}$. Then $\delta_n \rightarrow 0$ (since $u_n < u_{n+1} < 3u_n$ for $n \geq n_0$), and $$\eps_n = \Big( \sum_{j= \psi(n)}^{n-1} \eps_j \Big) + \delta_n u_{n+1}.$$

Let us consider an increasing sequence $(\uti_n )_{n \geq 1}$ of positive integers such that $\uti_{n+1} = 2 \uti_n - \uti_{\psi(n)}$ for any $n \geq n_0$. As recalled in Lemma \ref{lem62bis}, Equation  \eqref{eq623} holds for the sequence $(\uti_n)$. The associated sequence $(\deltati_n)$ is identically zero, and $(\epsti_n)$ satisfies $\epsti_n = \sum_{j= \psi(n)}^{n-1} \epsti_j $, for $n \geq n_0$.

We consider now the quotient $\al_n = \eps_n / \epsti_n$. We have 
$$\sum_{j= \psi(n)}^{n-1} \al_j \epsti_j = \sum_{j= \psi(n)}^{n-1} \eps_j = \eps_n - \delta_n u_{n+1}
= \al_n \epsti_n - \delta_n u_{n+1}$$
hence
\begineq \label{eq625}
\al_n = \sum_{j= \psi(n)}^{n-1} \frac{\epsti_j}{\epsti_{\psi(n)} + \ldots + \epsti_{n-1}} \al_j + \delta'_n \al_n
\eneq
by letting
$$\delta'_n = \delta_n \frac{u_{n+1}}{\eps_n} = \delta_n \frac{1}{1 - \frac{u_n}{u_{n+1}}}.$$
Let $0 < \vareta < 1$ be such that $u_{n+1} \geq \frac{1}{1-\vareta} u_n$ for any $n \geq n_0$. Then we have $|\delta'_n| \leq \vareta^{-1} |\delta_n| $ hence $\delta'_n \to 0$. Letting $\al'_n = (1 - \delta'_n) \al_n$, Equation \eqref{eq625} yields 
\begineq \label{eq626}
\al'_n = \sum_{j= \psi(n)}^{n-1} \frac{\epsti_j}{\epsti_{\psi(n)} + \ldots + \epsti_{n-1}} \al_j .
\eneq

\medskip

Since $\psi \in \calf$, there exists  $B$ be such that $\psi(n) \geq n-B$ for any $n \geq n_0$. For $n \geq n_0$, let $\Inmu$ denote the convex hull of $\al_{n-B}$, \ldots, $\al_{n-1}$ in $\R$, say $\Inmu = [\zeta_{n-1} \al_{n-1}, \zeta'_{n-1} \al_{n-1}]$ with $\zeta_{n-1} \leq 1 \leq \zeta'_{n-1}$. 
Let $\mu_n = \max(1 - \delta'_n, \frac{1}{1 - \delta'_n})$ (since we may assume $|\delta'_n| \leq 1/2$ for any $n \geq n_0$). Then we have $\mu_n \geq 1$ and $\mu_n \to 1$. 

Let us prove (by induction on $\ell$) that for any integers $n \geq n_0$ and $\ell \geq -1$ we have
\begin{equation} \label{eqet23}
\al_{n+\ell} \in \Big[ \frac{1}{\mu_n \mu_{n+1} \ldots \mu_{n+\ell} } \frac{1 + (B^{\ell+1}-1)\zeta_{n-1}}{B^{\ell+1}} \al_{n-1} , \, 
 \mu_n \mu_{n+1} \ldots \mu_{n+\ell} \frac{1 + (B^{\ell+1}-1)\zeta'_{n-1}}{B^{\ell+1}} \al_{n-1} \Big].
 \end{equation}
This is true for $\ell = -1$; let $\ell \geq 0$ be such that \eqref{eqet23} holds for any $\ell' \in \{-1, \ldots, \ell-1\}$ (with the same $n$, which remains fixed). For any $j$ such that $\psi(n+\ell) \leq j \leq n + \ell - 2$, we have 
$$\al_j \geq  \frac{\zeta_{n-1} \al_{n-1}}{\mu_n \mu_{n+1} \ldots \mu_{n + \ell - 1}}$$
by using the definition of $\Inmu$ if $j \leq n-1$, and the induction hypothesis otherwise. For $j = n + \ell -1$, the  induction hypothesis reads
$$\al_{n+\ell-1} \geq  \frac{1}{\mu_n \mu_{n+1} \ldots \mu_{n + \ell - 1}}
\frac{1 + (B^{\ell}-1)\zeta_{n-1}}{B^{\ell}} \al_{n-1}.
$$
Now Equation   \eqref{eq626} implies that $\al'_{n+\ell}$ is a weighted average of $\al_{\psi(n+\ell)}$, \ldots, $\al_{n+\ell-1}$ with non-negative weights; moreover the weight of $\al_{n+\ell-1}$ is $\epsti_{n+\ell-1} / (\epsti_{\psi(n+\ell)} + \ldots + \epsti_{n+\ell-1}) \geq 1/B$ since the sequence $(\epsti_n)_{n \geq n_0}$ is non-decreasing. Therefore the previous lower bounds yield:
$$\al'_{n+\ell} \geq \frac{1}{\mu_n \mu_{n+1} \ldots \mu_{n+\ell-1} } \frac{1 + (B^{\ell+1}-1)\zeta_{n-1}}{B^{\ell+1}} \al_{n-1} .$$
Now
$$\al_{n+\ell} = \frac{1}{1 - \delta'_{n+\ell}} \al'_{n+\ell} \geq \frac{1}{\mu_{n+\ell}} \al'_{n+\ell},$$
thereby proving the lower bound in \eqref{eqet23}. The upper bound can be proved in the same way, concluding the proof of  \eqref{eqet23}.

\medskip

Applying \eqref{eqet23}  for $\ell = 0, 1, \ldots, B-1$ yields
$$1 \leq \frac{\zeta'_{n+B-1}}{\zeta_{n+B-1}} \leq \mu_n^2 \mu_{n+1}^2 \ldots \mu_{n+B-1}^2 \frac{1 + (B^B-1)\zeta'_{n-1}}{1 + (B^B-1)\zeta_{n-1}}.$$
From this inequality we obtain
\begin{eqnarray*}
0 \leq   \frac{\zeta'_{n+B-1}}{\zeta_{n+B-1}} -1 
&\leq &  \mu_n^2 \mu_{n+1}^2 \ldots \mu_{n+B-1}^2 - 1 +  
\frac{\mu_n^2 \mu_{n+1}^2 \ldots \mu_{n+B-1}^2
(B^B-1)  (\zeta'_{n-1} - \zeta_{n-1})}{1 + (B^B-1)\zeta_{n-1}} \\
&\leq &  \mu_n^2 \mu_{n+1}^2 \ldots \mu_{n+B-1}^2 - 1 +  \mu_n^2 \mu_{n+1}^2 \ldots \mu_{n+B-1}^2
\frac{B^B-1}{B^B} \Big(  \frac{\zeta'_{n-1}}{\zeta_{n-1}} -1 \Big) .
\end{eqnarray*}
Since $\mu_n \to 1$, this proves that   $ \zeta'_{n+\ell B-1} / \zeta_{n+\ell B-1}$ tends to 1 as $\ell$ tends to infinity, from which we deduce  $\al_{n+\ell B-1} / \al_{n+\ell B-2} \rightarrow_{\ell \to \infty} 1$. Taking $B$ consecutive values for $n$ yields $\al_\ell / \al_{\ell-1} \rightarrow_{\ell \to \infty} 1$. Since $1 \leq \epsti_n / \epsti_{n-1} \leq B$ for any $n \geq n_0$, this implies $\eps_n / \eps_{n-1} = \epsti_n / \epsti_{n-1} + \odu$, hence
\begineq \label{eq630}
\frac{\eps_n }{ \eps_{n-t}} = \frac{\epsti_n}{ \epsti_{n-t}} + \odu \mbox{ for any fixed } t, \mbox{ as } n \to \infty.
\eneq
Moreover the quantities $\eps_n / \eps_{n-t}$ and $\epsti_n /  \epsti_{n-t}$ remain (as $n$ varies) in a fixed interval $[a_t, b_t] \subset (0, + \infty)$ depending only on $t$, so that we have also 
$\frac{\eps_{n-t} }{ \eps_{n}} = \frac{\epsti_{n-t}}{ \epsti_{n}} + \odu $. 

\medskip

For any $T \geq 1$, let us consider 
$$  \densi_T = \limsup_{n \to \infty} \frac{u_{n+1} - u_{n-T}}{u_n - u_{n-T}} \mbox{ and }
\widetilde \densi_T = \limsup_{n \to \infty} \frac{\uti_{n+1} - \uti_{n-T}}{\uti_n - \uti_{n-T}}.$$
The previous relations, and the definitions of $\eps_n$ and $\epsti_n$,  yield
$$\frac{u_{n+1} - u_{n-T}}{u_n - u_{n-T}} = 1 + \frac{1}{\frac{\eps_{n-T}}{\eps_n} + \ldots + \frac{\eps_{n-1}}{\eps_n}} = 
1 + \frac{1}{\frac{\epsti_{n-T}}{\epsti_n} + \ldots + \frac{\epsti_{n-1}}{\epsti_n} + \odu} 
=  \frac{\uti_{n+1} - \uti_{n-T}}{\uti_n - \uti_{n-T}} + \odu.$$
Letting $n$ tend to infinity, we obtain $ \densi_T  = \widetilde \densi_T$ for any $T \geq 1$. Now we are going to compute, separately, the limits of $ \densi_T$ and $ \widetilde \densi_T$ as $T$ tends to infinity.

Let us start with $ \densi_T$. We have
$$\frac{u_{n+1} - u_{n-T}}{u_n - u_{n-T}}  = \frac{  \frac{u_{n+1}}{u_n} -   \frac{u_{n-T}}{u_n}}{1 - 
   \frac{u_{n-T}}{u_n}}$$
   with $0 \leq \frac{u_{n-T}}{u_n} \leq (1 - \vareta)^T$ for any $n \geq n_0+T$, hence
   $$ \frac{u_{n+1}}{u_n} - (1 - \vareta)^T \leq \frac{u_{n+1} - u_{n-T}}{u_n - u_{n-T}}  \leq 
   \frac{1}{1 - (1 - \vareta)^T} \frac{u_{n+1}}{u_n} .$$
Lettting $n$ tend to infinity, we obtain:
$$\Big( \limsup_{n \to \infty}   \frac{u_{n+1}}{u_n} \Big) - (1 - \vareta)^T \leq \densi_T  \leq 
   \frac{1}{1 - (1 - \vareta)^T} \Big( \limsup_{n \to \infty}   \frac{u_{n+1}}{u_n} \Big),$$
   thereby proving that
$$ \lim_{T \to \infty} \densi_T =     \limsup_{n \to \infty}   \frac{u_{n+1}}{u_n} .$$

Now lemma \ref{lem62bis} yields $\liminf  \frac{\uti_{n+1}}{ \uti_n} > 1$, therefore the same arguments
give
$$ \lim_{T \to \infty} \widetilde \densi_T =     \limsup_{n \to \infty}   \frac{\uti_{n+1}}{\uti_n} = \densi(\psi) .$$

Since  $ \densi_T  = \widetilde \densi_T$ for any $T$, invoking the unicity of this limit enables us  to conclude the proof of    Lemma \ref{lem62}.

\subsection{The Palindromic Prefix Method} \label{subsec62}

Let $w$ be an infinite non ultimately periodic word such that $\densi(w) < 2$. As in \S \ref{subsec26}, we denote by $(n_i)_{i \geq 1}$ the increasing sequence of all lengths of palindromic prefixes of $w$. As mentioned in \S  \ref{subsec26}, the word $w$ can be written on a finite alphabet $\alphb$. 

\smallskip

We  consider an injective map $\varphi: \alphb \to \Netoile$   (which is just another way of saying that $w$ can be written on the alphabet of positive integers). Roy's idea \cite{RoyCRAS} is to consider the real number $\xi_{w, \varphi} $ defined by the infinite continued fraction expansion 
$$\xi_{w, \varphi} = [0, \varphi(w_1), \varphi(w_2), \varphi(w_3), \ldots] = 0 + \frac{1}{\varphi(w_1) + 
\frac{1}{\varphi(w_2) + \ldots}}.$$
Since $w$ is infinite and not ultimately periodic, $\xi_{w, \varphi} $ is neither rational nor quadratic.

Let $p_n / q_n$ denote the $n$-th convergent of $\xi_{w, \varphi} $; 
the classical properties of continued fraction expansions 
  (see for instance \cite{Schmidt}, Chapter I) give $\vert q_n \xi - p_n \vert \leq q_n ^{-1}$ and:
\begin{equation} \label{eq631}
\matri{q_n}{q_{n-1}}{p_n}{p_{n-1}} = 
\matri{\varphi(w_1)}{1}{1}{0}
\matri{\varphi(w_2)}{1}{1}{0} \ldots
\matri{\varphi(w_n)}{1}{1}{0}. 
\end{equation}
For $n=n_i$, the finite word $w_1 w_2 \ldots w_n$ is  a palindrome, so that this product is a symmetric matrix, hence $q_{n_i-1} = p_{n_i}$. In this case, we have simultaneous rational approximants to $\xi_{w, \varphi}$ and $\xi_{w, \varphi}^2$ with the same denominator, namely:
\begin{equation} \label{eq3624}
\left\{
\begin{array}{l}
|q_{n_i} \xi_{w, \varphi} - p_{n_i}| \leq q_{n_i}^{-1} \\
|q_{n_i} \xi_{w, \varphi}^2 - p_{n_i-1}| = |(q_{n_i} \xi_{w, \varphi} - p_{n_i}) \xi_{w, \varphi} + (q_{n_i-1} \xi_{w, \varphi} - p_{n_i-1})| \leq (1+\xi_{w, \varphi})q_{n_i-1}^{-1}.
\end{array}
\right.
\end{equation}
This gives $L(\vs_i) \leq c V_i^{-1}$ where $\vs_i = (q_{n_i}, p_{n_i}, p_{n_i-1}) \in \Z^3$ is the vector corresponding to the symmetric matrix \eqref{eq631}, and $c$ is a constant depending only on $w$ and $\varphi$. The definition of $\bezprefpal$  gives immediately (as in \cite{CRASasim}):
\begin{equation} \label{eq2525}
\bezprefpal \leq \limsup_{i \to \infty } \frac{\log V_{i+1}}{\log V_i}. 
\end{equation}
Our main result asserts that equality holds, and gives the value of this upper limit in  terms of $w$:

\begin{Th} \label{th63}
We have $ \bezprefpal =   \densi(w) $.
\end{Th}

The proof also shows that the function $\psi$ associated with $\xi$ (in Part $(a)$ of Theorem \ref{th51}) is the same, modulo $\calr$, as the one associated with $w$ (in Theorem \ref{thcombi}). 

\bigskip

\Question Does the equality $ \beuprefpal =   \densi(w) $ hold ? 

\medskip

When $w$ is the Fibonacci word (considered by Roy,  see Example \ref{explefibo}), Davenport-Schmidt's lower bound yields $\nbor \leq   \beuprefpal \leq  \bezprefpal = \densi(w) = \nbor$, hence equality holds. For any characteristic Sturmian word, Bugeaud and Laurent have proved  \cite{BL}   that equality always holds.  Using the methods introduced here, it is possible to generalize this result to any word $w$ equal modulo $\calsim$ to a characteristic Sturmian word (where $\calsim$ was defined in \S \ref{subsec26}). 

\bigskip

The proof of Theorem \ref{th63} falls into two parts. First, the analytic study (\S \ref{subsec61}) of the asymptotics of the recurrence relation associated with $w$ enables us to prove that $\limsup_{i \to \infty } \frac{\log V_{i+1}}{\log V_i} = \densi(w)$. Then Equation \eqref{eq2525} implies the upper bound in Theorem \ref{th63}: this upper bound is only of combinatorial, and analytic, nature. For the lower bound, one has to apply Theorem \ref{th51}: diophantine properties really come into the play.

\bigskip

In the proof of Theorem \ref{th63}, we shall use twice the following lemma (see \cite{SFcombi}, \S 5.1):

\begin{Lemme} \label{lem64}
Let $i_0$, $i_1$ and $i_2$ be positive integers such that $\min(n_{i_0}, n_{i_1}) \leq n_{i_2} \leq n_{i_0}+n_{i_1}$. Then the following statements are equivalent:
\begin{itemize}
\item The prefix of $w$ with length $n_{i_0} + n_{i_1} - n_{i_2}$ is a palindrome.
\item The vectors $\vs_{i_0}$, $\vs_{i_1}$ and $\vs_{i_2}$ are linearly dependent.
\end{itemize}
Moreover, 
\begin{itemize} 
\item If $i_2 \geq i_1 $ and $n_{i_2} \geq n_{i_0} - n_{i_1}$ then these statements hold.
\item If these statements hold then denoting by $i$ the integer such that $n_i = n_{i_0} + n_{i_1} - n_{i_2}$ we have 
\begineq \label{eq641}
[\vs_{i_0}, \vs_{i_1}, \vs_{i_2}] = \pm \, \vs_{i}
\eneq
and $\pal_i = \pal_{i_0} \pal_{i_2}^{-1} \pal_{i_1}$, where $\pal_j$ is the palindromic prefix of $w$ with length $n_j$.
\end{itemize}
\end{Lemme}

\Dem of Lemma \ref{lem64}:
To begin with, let us assume $i_2 \geq i_1 $ and $n_{i_2} \geq n_{i_0} - n_{i_1}$, and prove that the prefix of $w$ of length $n_{i_0} + n_{i_1} - n_{i_2}$ is a palindrome. Let $1 \leq p \leq (n_{i_0}+n_{i_1}-n_{i_2})/2$; then we have $p \leq n_{i_1}$ hence:
$$w_{n_{i_0}+n_{i_1}-n_{i_2}+1-p} = w_{n_{i_2}-n_{i_1}+p} = w_{n_{i_1}+1-p} = w_{p}$$
by using successively that $\pal_{i_0}$, $\pal_{i_2}$ and $\pal_{i_1}$ are palindromes.

Let us prove now that the first two statements are equivalent, and at the same time that \eqref{eq641} holds. We use the crucial formula (introduced in this context by Roy \cite{RoyPLMS}): 
$$\det(M_{i_0}, M_{i_1}, M_{i_2}) = \trace(J M_{i_0} J M_{i_2} J M_{i_1})$$
where $M_i$ is the symmetric matrix that corresponds to $\vs_i$, and $J={\tiny \left[ \begin{array}{cc} 0 & 1 \\ -1 & 0 \end{array} \right]}$. Now a matrix $M$ is symmetric if, and only if, $\trace(JM) = 0$. Therefore the linear dependence of $\vs_{i_0}$, $\vs_{i_1}$ and $\vs_{i_2}$ means that the matrix $M_{i_0} J M_{i_2} J M_{i_1}$, which is equal to $\pm M_{i_0} M_{i_2}^{-1} M_{i_1}$, is symmetric.

On the other hand, the prefix of $w$ with length $n_{i_0} + n_{i_1} - n_{i_2}$ is a palindrome if, and only if, the matrix $M$ given by Equation \eqref{eq631} with $n = n_{i_0} + n_{i_1} - n_{i_2}$ is symmetric (by unicity of such a decomposition). To conclude the proof, it is therefore enough to prove the following equality:
\begineq \label{eq642}
M = M_{i_0} M_{i_2}^{-1} M_{i_1}.
\eneq
To prove \eqref{eq642}, we may assume $i_1 \leq i_2$ (by interchanging $i_1$ and $i_0$ if necessary). 
Then $\pal_{i_1}$ is a prefix of $\pal_{i_2}$, hence there is a word $b$ such that $\pal_{i_2} = \pal_{i_1}b$. The word $b$ is a suffix of $\pal_{i_2}$, therefore its mirror image 
$\miroir{b}$ is a prefix of $\pal_{i_2}$ (hence of $w$) since $\pal_{i_2}$ is a palindrome (i.e., $\pal_{i_2} = \miroir{\pal_{i_2}}$). Now $\miroir{b}$ consists in $n_{i_2} - n_{i_1} \leq n_{i_0}$ letters, therefore $\miroir{b}$ is a prefix of $\pal_{i_0}$. As $\pal_{i_0}$ is a palindrome, $b$ is a suffix of $\pal_{i_0}$: there is a word $c$ such that $\pal_{i_0} = cb$. Finally, we have $\pal_{i_0} \pal_{i_2}^{-1} \pal_{i_1} = cb(\pal_{i_1}b)^{-1} \pal_{i_1} = c$. As $c$ is the prefix of $w$ of length $n= n_{i_0} + n_{i_1} - n_{i_2}$, Equation \eqref{eq642} is proved. This concludes the proof of Lemma \ref{lem64}.

\medskip

\bigskip

\Dem of Theorem \ref{th63}: 
We consider the sequence $(\vs_i)$ defined right after Equation \eqref{eq3624}, and we let $\psi \in \calf$  be a function associated with $w$ as in   Theorem \ref{thcombi}.   In particular, we have $\densi(\psi) = \densi(w)$. 

Let $i$ be sufficiently large. Since $n_{i+1} = 2 n_i - n_{\psi(i)}$, it  follows from  Lemma \ref{lem64} (applied with $i_0 = i_1 = i$ and $i_2 = i+1$) that $[\vs_i, \vs_i, \vs_{i+1}] = \pm \vs_{\psi(i)}$. 

In order to apply Proposition \ref{proplem61}, we need rough estimates for $V_i$. Let $\us$ be a product of $n$ matrices of the shape $\matri{a}{1}{1}{0}$ with $1 \leq a \leq \Omega$. Then it is not difficult (see for instance \cite{BL}, Equation (11) p. 788) to prove that $2^{[n/2]} \leq U \leq (2 \Omega)^n$. Letting
$\Omega = \max \varphi(\alphb)$, this implies $V_i \leq (2 \Omega)^{n_i}$ and $V_{i+1} \geq V_i \, 2^{[(n_{i+1}-n_i)/2]}$. Now Lemma \ref{lem62bis} yields $\liminf n_{i+1} / n_i > 1$, hence  $V_{i+1} > V_i^\delta$ for some $\delta > 1$. Therefore Proposition \ref{proplem61} applies, and gives
$$\limsup_{i \to \infty } \frac{\log V_{i+1}}{\log V_i} = \densi(\psi) =\densi(w).$$

Now the points $\vs_i$ are primitive  (since $\gcd(q_{n_i}, p_{n_i}) = 1$), have positive first coordinates $q_{n_i}$, and are ordered according to these first  coordinates.  Moreover, let $i$ be such that 
$\vs_{i-1}$, $\vs_i$ and $\vs_{i+1}$ are linearly dependent. Lemma  \ref{lem64}, applied with  $i_0 = i-1$, $i_1 = i+1$ and $i_2 = i$, proves that $n_{i+1} + n_{i-1} - n_i$ is the length of a palindromic prefix of $w$. Since this length is strictly between $n_{i-1}$ and $n_{i+1}$, it is equal to $n_i$, so that $n_{i+1} = 2 n_i - n_{i-1}$ and $\psi(i) = i-1$. Then, as noticed in the beginning of the proof, $[\vs_i, \vs_i, \vs_{i+1}]$ is collinear to $\vs_{i-1}$. Using Equation \eqref{eq211}, this is enough to apply Part $(b)$ of Theorem \ref{th51}, and conclude the proof of Theorem \ref{th63}.

\subsection{Proof of Theorem \ref{th1}} \label{subsec53}

Let $\xi$ be an irrational  non-quadratic real number such that $\bez < 2$. 
Denote by $(u_i)$ the sequence in Part $(a)$ of Theorem \ref{th51} (with $\eps = \eps_1/2$, say), and by  $\psi \in \calf$ the associated function. Using also  Proposition \ref{proplem61}, we have $\bez = \densi(\psi)$. Now, choosing a word $w$ associated with $\psi$ as in Theorem \ref{thcombi}, we have $\densi(w) = \densi(\psi)$ (see \S \ref{subsec26}), so that $\bez  = \densi(w)$.

Conversely, if $w$ is a word such that $1 < \densi(w) < 2$, then $w$ is not ultimately periodic and Roy's palindromic prefix method provides (for any embedding  $\varphi: \alphb \to \Netoile$) an irrational non-quadratic real number  $\xi_{w, \varphi} $ for which $\bezprefpal = \densi(w)$ thanks to Theorem \ref{th63}. 

This concludes the proof of Theorem \ref{th1}.

\section{Open Questions} \label{sec7}

Given a non ultimately periodic word $w$ such that $\densi(w) < 2$ (or an associated function $\psi \in \calf$),
it would be interesting to study more precisely the set of all numbers $\xi$  that correspond to $\psi$ (as in Theorem \ref{th51}). They all satisfy $\bez = \densi(w)$; but to  which extent do they all ``essentially'' come from the palindromic prefix method applied to  $w$? Answering this question would lead to a better understanding of the numbers $\xi$ such that $\bez < 2$. For instance, does there exist a $\xi$ such that $\beu < \bez < 2$ ? More generally, what can be said about the function $\eps \mapsto \beeps$ ? 

\smallskip

The methods used in this paper do not apply to numbers $\xi$ such that $\bez$ is finite but greater than (or equal to) 2. For instance, it is a completely open problem to determine the set of values $\geq 2$ taken by the exponent $\bez$.

\smallskip

The ``Fibonacci sequences'' introduced in \cite{Roytwoexp} provide, when $w$ is the Fibonacci word, an arithmetic enrichment of the palindromic prefix method. How does this enrichment generalize to other words $w$ such that $\densi(w) < 2$ ? Working this out would yield new examples of numbers $\xi$ such that $\beu < 2$, and would be a first step towards generalizing the results proved in this paper to the exponent $\beu$.

\smallskip

At last, Jarn\'{\i}k has proved \cite{Jarnik} (see also \cite{Roytwoexp} and \cite{BL}, end of \S 7) that $\beu$ is intimately connected to the ``dual'' problem of finding polynomials of degree 2, with not too large integer coefficients, that assume a small value at the point $\xi$. Maybe some exponents (depending on $\eps \in [0,1]$) can be defined and studied in relation with this dual problem.

\newcommand{\url}{\texttt}

\providecommand{\bysame}{\leavevmode ---\ }
\providecommand{\og}{``}
\providecommand{\fg}{''}
\providecommand{\smfandname}{\&}
\providecommand{\smfedsname}{\'eds.}
\providecommand{\smfedname}{\'ed.}
\providecommand{\smfmastersthesisname}{M\'emoire}
\providecommand{\smfphdthesisname}{Th\`ese}

\bigskip

\hspace{-\parindent}St\'ephane Fischler

\hspace{-\parindent}\'Equipe d'Arithm\'etique et de G\'eom\'etrie Alg\'ebrique

\hspace{-\parindent}B\^atiment 425

\hspace{-\parindent}Universit\'e Paris-Sud

\hspace{-\parindent}91405 Orsay Cedex, France

\hspace{-\parindent}stephane.fischler@\null {math}.u-psud.fr

\end{document}